\documentclass[final]{amsart}
\usepackage{amsfonts,amsmath,amssymb,mathtools}
\usepackage{graphicx}
\usepackage[comma, authoryear]{natbib}
\usepackage[colorlinks,citecolor=blue]{hyperref}
\usepackage{pdflscape,pdfsync}
\usepackage{xspace,colortbl,multirow}
\usepackage[dvipsnames]{xcolor}
\usepackage{color}
\usepackage{comment}

\usepackage[normalem]{ulem}
\usepackage{setspace}
\newtheorem{thm}{Theorem}[section]
\newtheorem{lem}[thm]{Lemma}
\newtheorem{prop}[thm]{Proposition}
\newtheorem{cor}[thm]{Corollary}
\theoremstyle{definition}
\newtheorem{defn}{Definition}[section]
\newtheorem{rem}{Remark}[section]
\newtheorem{exmp}{Example}[section]
\newtheorem{hyp}{Hypothesis}[section]

\numberwithin{equation}{section}

\def\disp{\displaystyle}

\def\seq#1{\left\{{#1}\right\}_{n\ge 1}}
\def\crochet#1{\langle #1 \rangle}
\DeclareMathOperator{\sign}{sign}

\newcommand{\dB}{\ensuremath{\mathbb{B}}}
\newcommand{\dC}{\ensuremath{\mathbb{C}}}

\newcommand{\dF}{\ensuremath{\mathbb{F}}}

\newcommand{\dH}{\ensuremath{\mathbb{H}}}
\newcommand{\I}{\ensuremath{\mathbb{I}}}
\newcommand{\dJ}{\ensuremath{\mathbb{J}}}

\newcommand{\dM}{\ensuremath{\mathbb{M}}}

\newcommand{\dR}{\ensuremath{\mathbb{R}}}

\newcommand{\cF}{\ensuremath{\mathcal{F}}}
\newcommand{\cG}{\ensuremath{\mathcal{G}}}

\newcommand{\cJ}{\ensuremath{\mathcal{J}}}

\newcommand{\cM}{\ensuremath{\mathcal{M}}}

\newcommand{\E}{{\rm E}}

\def\E{ \, {\rm E}}

\def\i1{ [-\infty,\infty]}
\def\ve{\varepsilon}
\newcommand{\esp}[1]{{\rm E}\left[{#1}\right]}

\def\sn{\sqrt{n}\; }

\def\ns{{\lfloor n s \rfloor}}
\def\nt{{\lfloor n t \rfloor}}

\def\ns{{\lfloor n s \rfloor}}
\def\nt{{\lfloor n t \rfloor}}

\usepackage[color, notref, notcite]{showkeys}
\usepackage{enumerate}
\definecolor{labelkey}{rgb}{0.6,0,1}

\newcommand{\BR}[1]{\textcolor{red}{#1}}
\newcommand{\JV}[1]{\textcolor{magenta}{#1}}

\begin{document}

\title[CLT for martingales-III: discontinuous compensators]
{Central limit theorems for martingales-III: discontinuous compensators}

\author{Bruno R\'{e}millard}
\address{Department of Statistics and Business Analytic, UAE University, and Department of Decision Sciences, HEC Montr\'{e}al}
\email{bruno.remillard@hec.ca}

\author{Jean Vaillancourt}
\address{Department of Decision Sciences, HEC Montr\'{e}al\\
3000, che\-min de la C\^{o}\-te-Sain\-te-Ca\-the\-ri\-ne,
Montr\'{e}al (Qu\'{e}\-bec), Canada H3T 2A7} \email{jean.vaillancourt@hec.ca}

\thanks{Partial funding in support of this work was provided by UAEU SRCL 2025 and the Natural
Sciences and Engineering Research Council of Canada}

\date{June 27, 2025}

\begin{abstract}
We propose a new weak convergence theorem for martingales, under gentler conditions than the
usual convergence in probability of the sequence of associated quadratic variations. Its proof requires
the combined use of Skorohod's $\cJ_1$-topology and $\cM_1$-topology on the space of c\`adl\`ag trajectories.
The emphasis is on those instances where the sequence of martingales or its limit is a mixture of stochastic processes
with discontinuities. Alternative conditions are set forth in the special cases of arrays of discrete time martingale differences and
martingale transforms. Examples of applications are provided, notably when the limiting process is a Brownian subordinator.
\end{abstract}


\maketitle

\section{Introduction}

This is the third part of a trilogy, with \citet{Remillard/Vaillancourt:2024a, Remillard/Vaillancourt:2025}.
The common goal in all three papers is to identify easily verifiable conditions under which a sequence $\{M_n\}$
of square integrable martingales, with (predictable) compensator $\langle M_n \rangle$, converges weakly
to another square integrable martingale $M$ (written $M_n \stackrel{\cJ_1}{\rightsquigarrow} M$ when using
the Skorohod $\cJ_1$-topology). It is well known that $\langle M_n \rangle\stackrel{\cJ_1}{\rightsquigarrow} \langle M \rangle$
alone does not suffice to yield $M_n\stackrel{\cJ_1}{\rightsquigarrow} M$, even when $\langle M \rangle$ has continuous trajectories.
Under the weaker requirement $\langle M_n \rangle \stackrel{f.d.d.}{\rightsquigarrow} \langle M \rangle$, the case where
$\langle M \rangle$ has continuous trajectories was treated in \citet{Remillard/Vaillancourt:2024a}; while weaker alternatives to the
conclusion $M_n \stackrel{\cJ_1}{\rightsquigarrow} M$ were explored in \citet{Remillard/Vaillancourt:2025} when $\langle M \rangle$
is no longer continuous. In the present paper, we propose conditions under which $M_n \stackrel{\cJ_1}{\rightsquigarrow} M$ ensues
from $\langle M_n \rangle\stackrel{\cJ_1}{\rightsquigarrow} \langle M \rangle$.
Section \ref{sec:m1_clt} states our two main results on functional CLTs for real-valued martingales.
Section \ref{sec:arrays_jumps} provides several illustrative examples. Section \ref{sec:MPU} presents a method
for lifting CLTs for non-stationary arrays with discontinuous limits, from those with continuous limits already extant in the literature.
The specific properties of inhomogeneous L\'evy processes needed here are assembled in
Appendix \ref{app:auxresultsLevy}. Some useful but technical properties of $\cM_1$-topology, as well as auxiliary ones on
weak convergence and most of the terminology, can be found in Appendix \ref{app:M1Topology}.
Tightness is handled in Appendices \ref{app:continuousinverse} and \ref{app:TightnessTheMainLemma}.
The proofs of the main results are gathered in Appendix \ref{app:main_results}.


\section{Main results}\label{sec:m1_clt}

Let $D= D[0,\infty)$ be the space of real-valued c\`adl\`ag trajectories (right continuous with left limits everywhere).
All processes considered here have their trajectories in $D$ and are adapted to a filtration $\dF = (\cF_t)_{t\ge 0}$
on a probability space $(\Omega,\cF,P)$ satisfying the usual conditions (notably, right continuity of the filtration and completeness).
Trajectories in $D$ are usually noted $x(t)$ but occasionally $x_t$.
Suppose that $M_n$ is a sequence of $D$-valued square integrable $\mathbb{F}$-martingales started at $M_n(0)=0$.
Its quadratic variation is denoted by $[M_n]$, its predictable compensator by $A_n:=\langle M_n \rangle$,
its largest jump by $J_T(M_n) := \sup_{s\in [0,T]}|\Delta M_n(s)|$ with $\Delta M_n(s) := M_n(s)-M_n(s-)$, and
the inverse process for its compensator $A_n$ by $\tau_n(s) := \inf\{t\ge0; A_n(t)>s\}$.
Define the rescaled $\mathbb{F}_{\tau_n}$-martingale $W_n := M_n\circ \tau_n$,
with compensator $\langle W_n\rangle = A_n\circ \tau_n$.
Except in the case where each $M_n$ is everywhere continuous, $W_n$ is not in general a Brownian motion.
Obtaining a CLT therefore requires building such a Brownian motion,
possibly on an enlargement of the stochastic basis $(\Omega,\cF,\dF, P)$ with $\dF = (\cF_t)_{t\ge 0}$.
(Such enlargements will be used systematically in this paper and are understood to affect some statements implicitly,
without further ado, for instance in some of the proofs). Since the goal is to obtain a weak convergence theorem for sequences of the form $x_n\circ y_n$, it is necessary to study when this is possible, under either of the $\cJ_1$-topology or the $\cM_1$-topology;
see Appendix \ref{app:M1Topology} for more details. We next recall the main definitions, before giving an example showing that one cannot expect convergence of composition, except in very special cases, when $D$ is equipped with the $\cM_1$-topology.

\subsection{Main topologies for convergence of  c\`adl\`ag processes  }\label{ssec:topo}

Throughout this section $\dB$ is an arbitrary separable Banach space $\dB$.
Some notation must be introduced ---any unexplained terminology can be found in \citet{Ethier/Kurtz:1986},
\citet{Whitt:2002}, \citet{Jacod/Shiryaev:2003} or \citet{Remillard/Vaillancourt:2024a, Remillard/Vaillancourt:2025}.

Let $D= D[0,\infty)$ be the space of $\dB$-valued c\`adl\`ag trajectories (right continuous with left limits everywhere), for some separable Banach space $\dB$, i.e., $\dB$ is a complete separable normed linear space with norm $\|\cdot\|$. $\dB$ will usually be $d$-dimensional Euclidean space $\dR^d$ in the examples; it will always be specified implicitly by the context at hand so the subscript $\dB$ is omitted. All processes considered here have their trajectories in $D$ and are adapted to a filtration $\dF = (\cF_t)_{t\ge 0}$ on a probability space
$(\Omega,\cF,P)$ satisfying the usual conditions (notably, right continuity of the filtration and completeness).
Trajectories in $D$ are usually noted $x(t)$ but occasionally $x_t$.

On $\dB^3$, set $\dC(x_1,x_2,x_3)=\|x_3-x_1\|$, $\dJ(x_1,x_2,x_3) = \|x_2-x_1\|\wedge \|x_2-x_3\|$, where $k\wedge\ell = \min(k,\ell)$
and let $\dM(x_1,x_2,x_3)$ be the minimum distance between
$x_2$ and the Banach space segment $[x_1,x_3]:=\{\lambda x_1+(1-\lambda)x_3\in \dB; \lambda\in[0,1]\}$. Then
$\dM(x_1,x_2,x_3) = 0 $ if $x_2 \in [x_1,x_3]$ and otherwise $\dM(x_1,x_2,x_3) \le \dJ(x_1,x_2,x_3)$.

For $\dH=\dC$, $\dH=\dJ$ or $\dH=\dM$, $T>0$, and $x \in D$, set
$$
\omega_\dH(x,\delta,T) = \sup_{0 \le t_1 < t_2 < t_3\le T, \; t_3-t_1<\delta}
\dH\{x(t_1),x(t_2),x(t_3)\}.
$$
Using the terminology in \citet[Theorem 3.2.2]{Skorohod:1956}, a sequence of processes $X_n$ converges weakly
under the $\mathcal{H}_1$-topology on $D$ to $X$,
denoted $X_n\stackrel{\mathcal{H}_1}{\rightsquigarrow} X$ when $\dH=\dJ$ or $\dH=\dM$, if and only if
\begin{enumerate}
\item[(i)]
the finite dimensional distributions of $X_n$ converge to those of
$X$, denoted $X_n \stackrel{f.d.d.}{\rightsquigarrow} X$, over a dense set of times in $[0,\infty)$ containing $0$;

\item[(ii)]
 for any $\epsilon>0$, and $T>0$,
\begin{equation}\label{eq:modulus}
\lim_{\delta\to 0}\limsup_{n\to\infty}P\{\omega_\dH(X_n,\delta,T) >
\epsilon\}=0.
\end{equation}
\end{enumerate}
The $\cM_1$-topology is weaker than the more commonly used $\cJ_1$-topology,
as shown in \citet{Skorohod:1956}. Similarly write $X_n\stackrel{\mathcal{C}}{\rightsquigarrow} X$ for weak convergence
under the $\mathcal{C}$-topology on $D$, that of uniform convergence over compact time sets, defined by
taking $\dH=\dC$. Remember that, when limit $X$ is continuous,
$X_n\stackrel{\mathcal{H}_1}{\rightsquigarrow} X$ and
$X_n\stackrel{\mathcal{C}}{\rightsquigarrow} X$ are equivalent for either $\dH=\dJ$ or $\dH=\dM$.
However, when we only know that processes $X_n$ are continuous,
$X_n\stackrel{\mathcal{H}_1}{\rightsquigarrow} X$ and
$X_n\stackrel{\mathcal{C}}{\rightsquigarrow} X$ are equivalent for $\dH=\dJ$
(and then $X$ is continuous as well) but this is not so for $\dH=\dM$.
One can observe $X_n\stackrel{\cM_1}{\rightsquigarrow} X$ with continuous processes $X_n$
converging to a discontinuous limit $X$. We shall see important examples of this phenomenon in the context of CLTs.

The $\cM_1$-topology on $D[0,1]$ was introduced in \citet{Skorohod:1956} and extended to $D$ in \citet{Whitt:1971}.
By \citet[Section 12.8]{Whitt:2002}, this topology is generated by a metric with respect to which $D$ is Polish (separable and complete).
Here we use the characterization of $\cM_1$-convergence due to \citet[Theorem 2.4.1]{Skorohod:1956}, by way of the special modulus defined above \eqref{eq:modulus}, equivalent variants of which are discussed at length in \citet[Section 12.5]{Whitt:2002}.

We are now in a position to present the first example of non-convergence for the composition, under the $\cM_1$-topology.
\begin{exmp}\label{ex:counterexample}
Let $y_n(t)=t/(2-1/n)$, $t\in [0,2-1/n]$, $y_n(t)=2+n(t-2)$, $t\in [2-1/n,2]$, and $y_n(t)=t$ for all $t\ge2$;
$y(t)=t/2$, $t\in [0,2)$ and $y(t)=t$ for all $t\ge2$; $x(t)=0$ for $t\in [0,1]\cup [2,\infty)$, $x(t)=1-|2t-3|$, $t\in [1,2]$.
Not only is $y^{-1}$ continuous, since $y$ is strictly increasing everywhere,
but every $y_n$ is everywhere Lipschitz and strictly increasing as well. However, even though
$y_n \stackrel{\cM_1}{\rightsquigarrow} y$ holds, $x\circ y_n \stackrel{\cM_1}{\not \rightsquigarrow} x\circ y\equiv 0$,
in spite of $x$ being a Lipschitz function.
In fact,  $x\circ y_n(t) = 0$ if $t\in [0,2-1/n]\cup [2,\infty)$, $x\circ y_n(t) = 2-2n(2-t)$, $t\in [2-1/n, 2-1/(2n)]$, and
$x\circ y_n(t) = 2n(2-t)$, $t\in [2-1/(2n),2]$, so $\dM\left\{x\circ y_n(2-1/n),x\circ y_n(2-1/(2n)),x\circ y_n(2) \right\} = \dM(0,1,0)=1$.
Hence, $\omega_{\dM}(x\circ y_n,\delta,T) = 1$, whenever $n> 1/\delta$ and $T\ge 2$.
\end{exmp}

The following result, first proposed to us by Bouchra R. Nasri, provides necessary and sufficient conditions under which
the composition map is $\cM_1$-continuous. The proof can be found in \citet[Proposition A.2]{Remillard/Vaillancourt:2025}.

\begin{prop}\label{prop:Whitt_composition}
Suppose that $(x_n, y_n) \stackrel{\cM_1}{\to} (x, y) $ in $D\left([0,\infty):\dR^k\right) \times  D_{\uparrow}$, where $x$ is continuous. Then
$x_n \circ y_n \stackrel{\cM_1}{\to}  x\circ  y$ in $D\left([0,\infty):\dR^k\right)$ if and only if $x$ is monotone on $[y(t-), y(t)]$,
for any $t\in {\rm Disc}(y)$.
\end{prop}
In Example \ref{ex:counterexample}, the non-convergence of $x\circ y_n$ to $x\circ y$ makes sense
since $x $ is not monotone on $[1,2] = [y(2-),y(2)]$. Proposition \ref{prop:Whitt_composition} also leads to
the following, minimal requirements.

\subsection{Hypotheses and main results}\label{ssec:hyp}

Recall that $W_n := M_n\circ \tau_n$, with compensator $\langle W_n\rangle = A_n\circ \tau_n$.

\begin{hyp}\label{hyp:An}
All of the following hold:
\begin{enumerate}
\item[(a)] $A_n(t) = \langle M_n\rangle_t \to\infty$ as $t\to\infty$ almost surely, for each fixed $n\ge1$;
\item[(b)] There is a $D$-valued process $A$ started at $A(0)=0$, such that \\
(i) $A_n\stackrel{\cJ_1}{\rightsquigarrow} A$; \\
(ii) $\disp \lim_{n\to\infty}E\left\{A_n(t)\right\}=E\left\{A(t)\right\} < \infty$ for all $t\ge 0$; \\
(iii) $A(t) \to\infty$ as $t\to\infty$ almost surely.
\item[(c)] $\disp \lim_{n\to\infty}\esp{A_n\circ \tau_n\circ A_n(t)-A_n(t)}=0$ for all $t\ge0$;
\item[(d)] $\disp \lim_{n\to\infty}\esp{\langle W_n\rangle_t}  =t$ for all $t\ge0$.
\item[(e)] Any of the following three conditions holds: \\
(i) $A$ is an inhomogeneous L\'evy process; \\
(ii) $\tau(t) = \inf\{s\ge0; A(s)>t\}$ is an inhomogeneous L\'evy process such that
$\tau_n \stackrel{f.d.d.}{\rightsquigarrow} \tau$ and $A_n(0)=0$ for every $n$; \\
(iii) $\tau(t) = \inf\{s\ge0; A(s)>t\}$ is an inhomogeneous L\'evy process and $\tau_n(0)=0$ for every $n$.
\end{enumerate}
\end{hyp}

\begin{rem}
Hypothesis \ref{hyp:An} collectively forms a strengthened version of both \citet[Hypothesis 2.1]{Remillard/Vaillancourt:2024a}, where the further assumption of continuity for the limiting compensator $A$ eases several aspects of the proofs, and
\citet[Hypotheses 3.1, 3.2 and 3.3 combined]{Remillard/Vaillancourt:2025}, where convergence is achieved in a much weaker sense than here.
The main change is Hypothesis \ref{hyp:An}b.i, strengthened from $A_n \stackrel{f.d.d.}{\rightsquigarrow} A$ (or equivalently
$A_n\stackrel{\cM_1}{\rightsquigarrow} A$ in the present context), in the light of Proposition \ref{prop:Whitt_composition}.
Note that $0\le\tau_n\circ A_n(t)-t=\sup\{h\ge0; A_n(t+h)=A_n(t)\}$, so $A_n\circ\tau_n\circ A_n(t)\ge A_n(t)$.
It then follows from Hypothesis \ref{hyp:An}a
that $\tau_n\circ A_n(t)-t <\infty$. Furthermore,
$$
E\left[ M_n(t) \times W_n \circ A_n(t) \right]= E\left[ M_n(t) \times M_n\circ \tau_n\circ A_n(t) \right] = E\left[\left\{ M_n(t)\right\}^2\right],
$$
so
\begin{eqnarray}
\label{eq:next_jump}
E\left\{ |W_n\circ A_n(t)-M_n(t)|^2 \right\} &=&  E\left\{ |M_n\circ \tau_n\circ A_n(t)-M_n(t)|^2 \right\}  \nonumber\\
 & = & E\left[ \left\{ W_n\circ A_n(t)\right\}^2 -\left\{M_n(t)\right\}^2 \right] \nonumber \\
 & = & E\left\{ \langle M_n\rangle_{\tau_n\circ A_n(t)} - \langle M_n\rangle_t \right\} \nonumber \\
 &=& \esp{A_n\circ \tau_n\circ A_n(t)-A_n(t)}.
\end{eqnarray}
Therefore, Hypothesis \ref{hyp:An}c states that, by analogy with the continuous case,
$M_n$ is almost a Brownian motion evaluated at its compensator. Furthermore, the first jump of $M_n$ after a fixed time $t$,
if there is one, occurs at the end of a flat spot for $A_n$,
and its second moment converges to $0$ as $n$ goes to $\infty$.
This technical requirement is automatically met if any of the following conditions is satisfied:
$M_n$ is continuous everywhere for all $n$; $\tau_n\circ A_n(t)=t$ for all $t\ge0$ and $n$;
$A_n$ is strictly increasing everywhere for all $n$, which implies the preceding condition;
$\tau_n$ is continuous everywhere for all $n$. Note that if either $A_n$ or $\tau_n$ is continuous everywhere,
then the other ($\tau_n$ or $A_n$, respectively) is strictly increasing everywhere, and vice versa.
Asymptotic rarefaction of jumps hypotheses (stronger versions of Hypothesis \ref{hyp:An}c
impacted by all the largest jumps at once) were introduced in \citet{Rebolledo:1980},
where they were shown to be less stringent than the often quoted Lindeberg condition.
These hypotheses are formulated in terms of convergence in probability and the reader should consult \citet{Rebolledo:1980} for details.
Nevertheless, Hypothesis \ref{hyp:An}c will sometime fail for sequences and arrays of discrete martingales,
especially when they are converging to limits with discontinuous trajectories ---
see the end of Section \ref{sec:m1_clt} for an alternative CLT in this context and Section \ref{sec:arrays_jumps}  for some examples.
Hypothesis \ref{hyp:An}c is adapted to sequences of processes with trajectories devoid of horizontal sections or flat spots,
such as stochastic integrals over integrands and integrators that never remain null over any interval of time.
\end{rem}

\begin{rem}
The mapping $A_n\mapsto\tau_n$ on $D$ is not continuous in the $\cJ_1$-topology,
not even when the application of the  mapping is restricted to uniformly convergent
sequences $A_n$ --- for a counterexample, see \citet[Example 13.6.1]{Whitt:2002}.
However, $A_n\mapsto\tau_n$ is continuous in the weaker $\cM_1$-topology,
under some additional conditions described in Proposition \ref{prop:inversecontinuity}.
This $\cM_1$-continuity can be instrumental in getting the full knowledge of the law of the limiting martingale,
through the use of either Hypothesis \ref{hyp:An}e.ii or \ref{hyp:An}e.iii.
\end{rem}

Hereafter, denote by $[M_n,B_n]$, the cross quadratic variation of two square integrable martingales, and
by $\langle M_n,B_n  \rangle$, their predictable cross quadratic variation.

To offer some perspective, we first state the following CLT in the $\mathcal{C}$-topology
\citep[Theorem 2.1]{Remillard/Vaillancourt:2024a}, where $A$ is assumed to be continuous.

\begin{thm}\label{thm:jumps_vanish}
Assume that $A$ has continuous trajectories and  $A_n \stackrel{f.d.d.}{\rightsquigarrow} A$ holds  (instead of Hypothesis \ref{hyp:An}b.i). Assume that Hypotheses \ref{hyp:An}a, \ref{hyp:An}b.ii and \ref{hyp:An}b.iii also hold; that
$J_t(M_n) \stackrel{Law}{\rightsquigarrow} 0$ for any $t>0$; that there exists an
$\dF$-adapted sequence of $D$-valued square integrable martingales $B_n$ started at $B_n(0)=0$ so that
\begin{enumerate}
\item $(B_n, A_n) \stackrel{\mathcal{C}}{\rightsquigarrow} (B,A)$ holds, where $B$ is a Brownian motion
with respect to its natural filtration $\dF_B = \{\cF_{B,t}: \; t\ge 0\}$ and $A$ is $\dF_B$-measurable;
\item  $\langle M_n, B_n\rangle_t \stackrel{Law}{\rightsquigarrow} 0$, for any $t\ge 0$.
\end{enumerate}
Then $(M_n, A_n,B_n) \stackrel{\mathcal{C}}{\rightsquigarrow} (M,A,B)$ holds,
where $M$ is a continuous square integrable ${\cF}_{t}$-martingale with respect to (enlarged) filtration
with predictable quadratic variation process $A$.
Moreover, $M = W\circ A$ holds, with $W$ a standard Brownian motion which is independent of $B$ and $A$.
\end{thm}

We now state our main result, the proof of which is given in Appendix \ref{pf:thmmain_terr1}.
Important applications of this result are the construction of processes with L\'evy subordinator to Brownian motion, a class of processes used in financial modelling. See, e.g., Example \ref{exmp:linnik}.

\begin{thm}\label{thm:main_terr1}
Assume that Hypothesis \ref{hyp:An} holds. Then there holds $W_n \stackrel{\mathcal{C}}{\rightsquigarrow} W$ and
$(M_n,A_n,W_n) \stackrel{\cJ_1}{\rightsquigarrow} (M,A,W)$,
with $M = W\circ A$ and $W$ is a Brownian motion, independent of process $A$.
Under either Hypothesis \ref{hyp:An}e.ii or \ref{hyp:An}e.iii, there also follows
$(A_n,\tau_n,W_n) \stackrel{\cM_1}{\rightsquigarrow} (A,\tau,W)$ and $W$ is independent of the pair $(A,\tau)$.
If, in addition, there holds $J_t(M_n) \stackrel{Law}{\rightsquigarrow} 0$ for any $t>0$,
then so do both $(M_n, A_n, W_n) \stackrel{\mathcal{C}}{\rightsquigarrow} (M,A,W)$ and
$(M_n,A_n,\tau_n,W_n) \stackrel{\cM_1}{\rightsquigarrow} (M,A,\tau,W)$, with $(M,A,W)$ is continuous.
\end{thm}

Theorems \ref{thm:jumps_vanish} and \ref{thm:main_terr1} are not corollaries of one another.
Neither process $A$ nor its inverse $\tau$ need be L\'evy processes in Theorem \ref{thm:jumps_vanish},
so it does not ensue from Theorem \ref{thm:main_terr1}. Meanwhile, a discontinuous $A$ makes the sole
Theorem \ref{thm:jumps_vanish} insufficient in proving Theorem \ref{thm:main_terr1}.

In the literature, one often encounters statements about the finite dimensional convergence of the normalized process $M_n/\sqrt{A_n}$.
Its behaviour can be deduced from our Theorem \ref{thm:main_terr1}, and the proof of the following corollary is given in Appendix \ref{pf:cor:standardization}.

\begin{cor}\label{cor:standardization}
Under Hypothesis \ref{hyp:An}, $(A_n,M_n/\sqrt{A_n}) \stackrel{f.d.d.}{\rightsquigarrow} \left(A,M/\sqrt{A}\right)$
holds with the choice $0/0=0$ and $\left(A,M/\sqrt{A}\right) \stackrel{Law}{=} (A,W)$.
Provided $A$ is continuous, $\left(A_n,M_n/\sqrt{A_n}\right) \stackrel{\mathcal{C}}{\rightsquigarrow} (A,W)$ also holds.
\end{cor}

The failure of either Hypothesis \ref{hyp:An}c or \ref{hyp:An}d precludes the use of Theorem \ref{thm:main_terr1}.
For instance, this is the case for Hypothesis \ref{hyp:An}c when the eventual limit $M$ has an integer-valued jump component,
as do all queueing processes.
Hypotheses \ref{hyp:An}c and \ref{hyp:An}d can sometimes be bypassed
through the following alternative (when $A$ is strictly increasing but may be discontinuous).

\begin{defn}\label{def:increasing}
A process $V$ is said to be strictly increasing when $P(V_t-V_s>0)=1$ holds for all choices of $0\le s<t<\infty$.
This property typifies L\'evy subordinators amongst inhomogeneous L\'evy processes.
\end{defn}

\begin{cor}\label{cor:main_terr1b}
Assume that Hypotheses \ref{hyp:An}a, \ref{hyp:An}b and \ref{hyp:An}e.i hold; that $A$ is strictly increasing;
 that $\esp{A\circ\tau(t)}=t$ for all $t\ge0$; that $\langle W_n\rangle_t$ are uniformly integrable for each fixed $t>0$.
Then $(M_n,A_n) \stackrel{\cJ_1}{\rightsquigarrow} (M,A)$ holds,
with $M = W\circ A$ and $W$ is a Brownian motion, independent of process $A$.
Furthermore, the conclusions of Corollary \ref{cor:standardization} still hold.
\end{cor}
The proof is given is Appendix \ref{pf:cor:main_terr1b}.

\begin{rem}\label{rem:strict}
When $A$ is strictly increasing, $\tau$ is continuous; under either Hypothesis \ref{hyp:An}e.ii or \ref{hyp:An}e.iii,
non-decreasing process $\tau$ is also an inhomogeneous L\'evy process and therefore both $\tau$ and $A$ are deterministic.
Hence the exclusion of both Hypotheseses \ref{hyp:An}e.ii and \ref{hyp:An}e.iii from the statement of
Corollary \ref{cor:main_terr1b} is explicit, since instances of asymptotics with $A$ deterministic can be handled directly by other means.
\end{rem}

\begin{rem}
Note that $(W_n,\tau_n)$ is an $\mathbb{F}_{\tau_n}$-semimartingale
since both $W_n$ and $\tau_n$ are $\mathbb{F}_{\tau_n}$-adapted, $\tau_n$ has finite variation and
$W_n$ is an $\mathbb{F}_{\tau_n}$-martingale. Under Hypothesis \ref{hyp:An}, the limit $(W,\tau)$ is both
an $\mathbb{F}_{\tau}$-semimartingale and an inhomogeneous L\'evy process.
As documented in \citet{Jacod/Shiryaev:2003}, previous approaches to the CLT exploited these features
of the limit $(W,\tau)$ in order to achieve convergence, by way of the associated characteristic functions.
The intuitive conditions in both Theorem \ref{thm:main_terr1} and Corollary \ref{cor:main_terr1b} are often easier to check for many
applications.
\end{rem}


\section{Discrete martingale arrays and pure jump martingales}\label{sec:arrays_jumps}

In this section, we consider the simplest sequences of real-valued discontinuous square integrable $\mathbb{F}$-martingales $M_n$,
namely those with step-valued trajectories arising from the weak invariance principle for arrays of discrete time martingale differences.
For each $n\ge1$, let the $\{X_{n,j}\}_{j\ge 1}$ be square integrable martingale differences with respect to
a filtration such that $\cG_{n,k} \supseteq \sigma\{X_{n,j}; j\le k\}$.
Define $\disp M_{n}(t) := \sum_{j=1}^{\nt} X_{n,j}$, with starting value $M_{n}(0)=0$, so
$\disp [M_n]_t = \sum_{j=1}^{\nt}X_{n,j}^2$ and $\disp A_n(t) = \langle M_n \rangle _t = V_{n,\nt}$, with $[M_n]_0=A_n(0)=0$,
where $\disp V_{n,k} = \sum_{j=1}^{k}E\left(X_{n,j}^2|\cG_{n,j-1}\right)$.
By construction $(M_n,[M_n],A_n)$ is $D$-valued and we only need to complete the filtration $(\cG_{n,k})_{k\ge 1}$ in order
to ensure that the set-up corresponds to our requirements. In Section \ref{ssec:hist}, we recall so early work on a CLT for a sequence (Proposition \ref{prop:eagleson}) and we improve the result by proving a functional CLT. Next, in Section \ref{ssec:hyp2}, we state new hypotheses well-suited for partial-sum processes, and in Section \ref{ssec:linnik_gamma}, we show some applications with discontinuous limits.

\subsection{Early work}\label{ssec:hist}

In the special case where $A(t)=t$ everywhere and all random variables are built on a common probability space,
the usual requirements for a martingale CLT to hold
($M_n \stackrel{\mathcal{C}}{\rightsquigarrow} W$ with Brownian motion $W$ as the limit)
are those in \cite{McLeish:1974} and \cite[Theorem 3]{Aldous:1978}, namely, $A_n(t) \stackrel{Pr}{\rightsquigarrow} t$,
plus a uniform boundedness condition on the second moment such as
$\disp \sup_n E\left(\max_{j\le \nt}X_{n,j}^2\right)<\infty$, in both cases for each $t\ge0$.
CLTs involving explicit mixtures of laws date back to the same period --- here is a sample result \cite[Corollary 3.1]{Hall/Heyde:1980}.

\begin{prop}\label{prop:eagleson}
Assume that the $\sigma$-algebras are nested: $\cG_{n,k}\subseteq \cG_{n+1,k}$ for all $n\ge1$ and all $1\le k \le n$.
Assume that
$A_n(1)\stackrel{Pr}{\rightsquigarrow}\eta$ as $n\to\infty$, for some random variable $\eta$ such that  $P\{\eta\in(0,\infty)\}=1$.
Suppose also that, for any $\epsilon>0$,
$\disp \sum_{j=1}^{n} E\left\{X_{n,j}^2\I(|X_{n,j}|>\epsilon) |\cG_{n,j-1}\right\} \stackrel{Pr}{\rightsquigarrow} 0$, as $n\to\infty$.
Then $M_{n}(1)\stackrel{Law}{\rightsquigarrow}Z\sqrt{\eta}$, where $Z\sim N(0,1)$ is independent of $\eta$.
\end{prop}

In \citet[Corollary 3.1]{Hall/Heyde:1980} the conclusion is stated for (the stronger) stable convergence.
Here is an example where $\eta$ is non deterministic, inspired by \citet[Section 3.2, Remarks]{Hall/Heyde:1980}.

\begin{exmp}\label{exmp:non_determ}
Set $\{Y_i\}$ and $\{\psi_i\}$ iid with $\disp P(Y_i=\pm 1) = P(\psi_i=\pm 1) = \frac{1}{2}$,
the two sequences being independent of each other. Let $\disp B_n(t):=\frac{1}{\sn}\sum_{i=1}^\nt Y_i\psi_{i}$,
$\disp M_n(t):=\frac{1}{\sn}\sum_{i=1}^\nt Y_i Z_{i-1}$, where $\disp Z_i:=\left(\sum_{j=1}^{i} \frac{Y_j}{j}\right)$ and $Z_0:=1$.
For every $n\ge1$, write $\cG_{n,k} := \cF_k := \sigma\{Y_i,\psi_i; i\le k\}$ and
$\disp Z_\infty = \sum_{i=1}^\infty \frac{Y_i}{i}\in\cF=\sigma\{Y_i; i\ge 1\}$. It follows that
$\langle M_n,B_n  \rangle = [M_n,B_n]_t = 0$ and
$\disp \langle M_n \rangle_t = [M_n]_t = \frac{1}{n}\sum_{i=1}^\nt Z_{i-1}^2\to A(t) = t Z_\infty^2$, with $\disp E\left(Z_\infty^2\right)= \frac{\pi^2}{6}$.
As $|Z_i|\le 1+\log{i}$, $\disp J_t(M_n) = \max_{1\le i \le \nt}|Y_i Z_{i-1}|/\sn\le \left(1+\log{n}\right)/\sn \to 0$ as $n\to\infty$,
Theorem \ref{thm:jumps_vanish} yields $M_n \stackrel{\mathcal{C}}{\rightsquigarrow} W\circ A$.
Proposition \ref{prop:eagleson} also holds, as the Lindeberg-type condition ensues from
$\disp \frac{1}{n}\sum_{j=1}^{n} (\log j)^2\I(\log j >\epsilon\sqrt{n})=0$ for all $n$ large enough. Note that our result is stronger than the result obtained by Proposition \ref{prop:eagleson}.
\end{exmp}

\begin{rem}\label{rem:meas_constraint}
Some measurability constraint on $\eta$ or on the family of $\sigma$-algebras is necessary, otherwise the conclusion fails ---
see \cite[Remark 6.1]{Dvoretzky:1972} and \cite[Section 3.3, Example 4]{Hall/Heyde:1980} for counterexamples.
These counterexamples all involve a limit $A$ stochastically dependent on all of $\{X_{n,j}\}$, as does Example \ref{exmp:non_determ}.
The construction of a Brownian motion $B$ rendering the limit $A$ to be $\dF_B$-measurable, with natural filtration
$\dF_B = \{\cF_{B,t}: \; t\ge 0\}$, while keeping $A$ independent of $W$, frees our approach from such constraints
on the martingale sequence itself. Focusing instead on the structure of the limit, expands the usefulness and applicability of the CLT.
\end{rem}

\subsection{Hypotheses for pure jump arrays}\label{ssec:hyp2}

Getting back to cases where $A$ is discontinuous, Corollary \ref{cor:main_terr1b} offers complementary versions with discontinuous limits,
of earlier CLTs with continuous limits --- see \citet{Merlevede/Peligrad/Utev:2019} for a recent survey and Section \ref{sec:MPU} for some such versions. In the current discrete context, the following Hypothesis \ref{hyp:An2} implies Hypothesis \ref{hyp:An}.
\\

 \begin{hyp}\label{hyp:An2}
All of the following hold:
\begin{enumerate}
\item[(a)] $V_{n,\infty}=\infty$ almost surely for each fixed $n\ge1$;
\item[(b)] There is a $D$-valued process $A$ started at $A(0)=0$, such that \\
(i) $A_n\stackrel{\cJ_1}{\rightsquigarrow} A$; \\
(ii) $\disp \lim_{n\to\infty}E\left\{A_n(t)\right\}=E\left\{A(t)\right\} < \infty$ for all $t\ge 0$; \\
(iii) $A(t) \to\infty$ as $t\to\infty$ almost surely.
\item[(c)] $\disp \lim_{n\to\infty} E\left\{X_{n,\nt+1}^2\right\}=0$, for any $t\ge 0$;
\item[(d)] $E\left\{V_{n,V_n^{-1}(t)}\right\}\to t$ as $n\to\infty$, for any $t\ge 0$.
\item[(e)] Either of the following two conditions holds: \\
(i) $A$ is an inhomogeneous L\'evy process; \\
(ii) $\tau(t) = \inf\{s\ge0; A(s)>t\}$ is an inhomogeneous L\'evy process such that
$\tau_n \stackrel{f.d.d.}{\rightsquigarrow} \tau$.
\end{enumerate}
\end{hyp}

Quick verifications confirm that: Hypothesis \ref{hyp:An2}a implies Hypothesis \ref{hyp:An}a;
\ref{hyp:An2}b is the same as \ref{hyp:An}b; \ref{hyp:An2}c is equivalent to \ref{hyp:An}c;
\ref{hyp:An2}d is equivalent to \ref{hyp:An}d; \ref{hyp:An2}e.i is the same as \ref{hyp:An}e.i;
and \ref{hyp:An2}e.ii is the same as \ref{hyp:An}e.ii, as $A_n(0)=0$.
Finally, \ref{hyp:An}e.iii always fails for discrete martingale arrays, since $A_n(0)=0$ and $V_{n,\infty}>0$ together
imply $\tau_n(0)>0$ for all $n\ge1$, so the inverse mapping $A_n\mapsto\tau_n$
may not be continuous at $A_n$ in the $\cM_1$-topology.
This is the case here even though the mapping $\tau_n\mapsto A_n$ is continuous in the $\cM_1$-topology.
Note also that Hypothesis \ref{hyp:An2}d holds if for any sequence of stopping times
$\gamma_n\to\infty$, $\lim_{n\to\infty} E\left(X_{n,\gamma_n}^2\right)=0$,
because of
$$
V_{n,V_n^{-1}(t)} = A_n\circ \tau_n(t) \ge t
\ge V_{n,V_n^{-1}(t)-1} = V_{n,V_n^{-1}(t)}-\E\left(X_{n,V_n^{-1}(t)}^2|\cF_{n,V_n^{-1}(t)-1}\right),
$$
where $V_n^{-1}(t)=\inf\{k\ge 0; V_{n,k}> t\}$.
Here $\disp \tau_n\circ A_n(t) = \frac{1+\nt}{n}\in\left(t,t+\frac{1}{n}\right]$.

\begin{rem}
In the frequent instance where $X_{n,j}^2>0$ holds almost surely for all $j$ and $n$, there ensues
$J_T(\tau_n) \le 1/n$; under such circumstances, $\tau_n\stackrel{\cJ_1}{\rightsquigarrow} \tau$ implies
$\tau_n \stackrel{\mathcal{C}}{\rightsquigarrow} \tau$, so $\tau$ is continuous and $A$ strictly increasing everywhere.
Thus strengthening $\tau_n \stackrel{f.d.d.}{\rightsquigarrow} \tau$ to $\tau_n\stackrel{\cJ_1}{\rightsquigarrow} \tau$
in Hypothesis \ref{hyp:An2}e.ii invariably restricts $\tau$ to be a continuous non-decreasing L\'evy process and therefore deterministic.
Hence so is $A$. The more interesting limits here arise from instances when $\tau_n\stackrel{\cJ_1}{\rightsquigarrow} \tau$ fails.
\end{rem}

\begin{rem}
Proving Hypothesis \ref{hyp:An2}b or \ref{hyp:An}b is generally the heart of the matter and the hardest part of the verifications.
One can check the $\cJ_1$-tightness using Theorem \ref{thm:M1tightness} and, separately, prove
$A_n \stackrel{f.d.d.}{\rightsquigarrow} A$ using the Cram{\'e}r-Wold device,
which ensures the sufficiency of the convergence of linear combinations
\begin{equation}\label{cramer_wold}
\sum_{j=1}^\ell \omega_j \{A_n(s_j)-A_n(s_{j-1})\} \stackrel{Law}{\rightsquigarrow} \sum_{j=1}^\ell \omega_j \{A(s_j)-A(s_{j-1})\},
\end{equation}
for every $\ell\ge1$, all $(\omega_j)_{1\le j\le \ell}\in\dR^{\ell}$ and each choice of $0=s_0 \le s_1 \le \cdots \le s_{\ell} \le t$.
The same argument works for $\tau_n \stackrel{f.d.d.}{\rightsquigarrow} \tau$.
\end{rem}

To illustrate the choices made through the above hypotheses, consider their effect on a single sequence of real-valued discrete time martingales.

\begin{exmp}\label{exmp:discrete}
Set $\disp M_n(t) = \frac{1}{a_n}\sum_{k=1}^\nt \xi_k$, where $\xi_k$ is a martingale difference sequence
with $E\left(\xi_k^2\right)<\infty$ and $a_n\to\infty$. Then $M_n$ is a martingale with respect to $\cF_{n,t}=\cF_{\nt}$,
$\disp \langle M_n \rangle_t = A_n(t) = \frac{1}{a_n^2}\sum_{k=1}^\nt E(\xi_k^2|\cF_{k-1})$ and
$\disp [M_{n}]_t = \frac{1}{a_n^2}\sum_{k=1}^\nt \xi_k^2$ respectively.
Suppose that $A_n(1) \stackrel{Law}{\rightsquigarrow} \eta$
and $a_{\nt}^2/a_n^2 \to \ell(t)$, with $\ell$ continuous, unbounded and started at $\ell(0)=0$.
Then $A_n \stackrel{f.d.d.}{\rightsquigarrow} A$ with $A(t) = \ell(t)\eta$, so Hypotheses \ref{hyp:An2}a and \ref{hyp:An2}b hold,
provided $\eta>0$ almost surely and $E(\eta)<\infty$ both hold.
Note that $\tau_n(0)=1/n$ and $\tau_n\circ A_n(t) = (\nt+1)/n>t$ for all $t\ge0$, so Hypothesis \ref{hyp:An2}c becomes
$\lim_{n\to\infty}\frac{1}{a_n^2} E(\xi_n^2)=0$, because of $a_{\nt}^2/a_n^2 \to \ell(t)$,
provided we assume $E(\xi_k^2|\cF_{k-1})>0$ to avoid trivialities.
Next, observe that $A_n\circ \tau_n(t) = A_n\left(\lceil nt \rceil/n\right)$.
Therefore, Hypothesis \ref{hyp:An2}d is equivalent to $\disp A(t) = \frac{t\eta}{E(\eta)}$, which implies $\ell(t)E(\eta)=t$.
This suggests Hypothesis \ref{hyp:An2}d is only adapted to the single discrete martingale context when
$A_n \stackrel{\mathcal{C}}{\rightsquigarrow} A$ holds for $A$ a straight line with an $\cF_0$-measurable slope.
These hypotheses are to be compared with the classical ones for CLTs:
\begin{itemize}
\item[(L)] (Lindeberg-type condition)
\begin{equation}\label{eq:condL}
\frac{1}{a_n^2}\sum_{k=1}^n E\left\{\xi_k^2\I(|\xi_k|> a_n\ve)|\cF_{k-1}\right\} \stackrel{Pr}{\rightsquigarrow} 0 \quad \text{ for any }\ve>0;
\end{equation}
\item[(M)] (McLeish-type condition)
\begin{equation}\label{eq:condM}
\sup_{s\le t}|[M_{n}]_s-A_n(s)| \stackrel{Pr}{\rightsquigarrow} 0 \quad \text{ for any }\ve>0.
\end{equation}
\end{itemize}
Keep in mind that \eqref{eq:condL} implies \eqref{eq:condM} \citep[Theorem 2.23]{Hall/Heyde:1980}
when $A_n(t)$ is assumed to be tight for every $t>0$. More on these and related conditions can be found
in \citet{Hall/Heyde:1980}. We provide such a comparison next, in a simpler but still illustrative context.
\end{exmp}

\begin{exmp}\label{exmp:lindeberg}
Pursuant to the notation of Example \ref{exmp:discrete}, the $\{\xi_k\}$ are now independent random variables
such that, for each $k\ge1$, $P(\xi_k=\pm k^{\alpha/2})=1/(2k^{\beta})$ while $P(\xi_k=0)=1-1/k^{\beta}$,
with $\beta\ge0$ and $\delta = \alpha-\beta+1\ge 0$. Set $a_n^2 = n^{\delta}$ if $\delta>0$ and $a_n^2 = \log{n}$ if $\delta=0$.
When $\alpha\in[-1,0]$, $\sup_{k\ge1}|\xi_k|\le1$ holds and condition \eqref{eq:condL} in Example \ref{exmp:discrete} is satisfied;
hence so is condition \eqref{eq:condM}, no matter what the values of $\delta\in[0,1]$ and $\beta\in[0,1]$ are. We proceed with $\alpha>0$.
When $\delta=0$, condition \eqref{eq:condL} now reads
 $\disp (\log n)^{-1}\sum_{k=\lfloor (\log n)^{1/\alpha}\ve^{2/\alpha} \rfloor+1}^n k^{-1}\to 0$ for any $\ve >0$;
 when $\delta>0$, condition \eqref{eq:condL} reads instead
 $n^{-\delta}\sum_{k=\lfloor n^{\delta/\alpha}\ve^{2/\alpha} \rfloor+1}^n k^{\delta-1}\to 0$ for any $\ve >0$. Hence,
when $\alpha>0$, condition \eqref{eq:condL}  holds if and only if both $\beta<1$ and $\delta>0$.
When $\delta=0$, Hypothesis \ref{hyp:An2}.b.iii also fails, as $A_n(t) \to 1$ for all $t>0$ in this case;
but when $\delta>0$, Hypotheses \ref{hyp:An2}a and \ref{hyp:An2}b hold since it means that all of the following must hold:
using $\disp A_n(t)=n^{-\delta}\sum_{k=1}^\nt k^{\delta-1}$,
$A_n(t)\to\infty$ as $t\to\infty$, for each fixed $n\ge1$, plus
$A_n(t) \to A(t) = \dfrac{t^\delta}{\delta} $ as $n\to\infty$, for all $t\ge0$,
and $A(t)\to\infty$ as $t\to\infty$, which only holds provided $\delta>0$ as well.
Hypotheses  \ref{hyp:An2}a and \ref{hyp:An2}b are equivalent to $\delta>0$.
The behavior of $\disp [M_{n}]_t  = n^{-\delta}\sum_{k=1}^\nt k^{\alpha}\I(\xi_k\neq0)$ is markedly different from that of $A_n(t)=E[M_{n}]_t$
--- for instance, $[M_{n}]_t \stackrel{a.s.}{\rightsquigarrow} 0$, for each $t\ge0$, whenever $\beta>1$,
by an application of the Borel-Cantelli lemma to the sequence $\{\I(\xi_k\neq0)\}$.
When $\delta>0$, $Z_n(t) = [M_{n}]_t-A_n(t) = n^{-\delta}\sum_{k=1}^\nt k^{\alpha}\left\{\I(\xi_k\neq0)-k^{-\beta}\right\}$
is a martingale for each fixed $n$, and
$$
\langle Z_n \rangle_t  = {\rm Var}\{Z_n(t)\} = \frac{1}{n^{2\delta}}\sum_{k=1}^\nt k^{2\alpha}\left(\frac{1}{k^\beta}-\frac{1}{k^{2\beta}}\right),
$$
satisfies, for each fixed $t>0$, $\langle Z_n \rangle_t \to t^{2\alpha}/(2\alpha)$ if $\beta=1$, $\langle Z_n \rangle_t \to 0$ if $\beta<1$
and $\langle Z_n \rangle_t \to \infty$ if $\beta>1$; while $Z_{n}(t) \stackrel{a.s.}{\rightsquigarrow} -t^\delta/\delta$ when $\beta>1$.
Condition \eqref{eq:condM} implies $Z_n(t) \stackrel{Pr}{\rightsquigarrow} 0$ for each $t>0$, which fails as soon as $\beta>1$ and $\delta\ge0$;
while $\beta<1$ and $\delta>0$ together ensure $\langle Z_n \rangle_t \to 0$ for each $t>0$ and therefore condition \eqref{eq:condM} is valid as well,
by way of an application of Lemma \ref{lem:lenglart}.
Meanwhile, the conditions of the CLT in Theorem \ref{thm:jumps_vanish} boil down to $\delta>0$ together with
$\max_{1 \le k\le n}\xi_k^2/a_n^2 \stackrel{Law}{\rightsquigarrow} 0$.
Invoking Theorem \ref{thm:M1tightness}, they imply that $[M_n]$ is $\mathcal{C}$-tight as well;
hence, $\langle Z_n \rangle_t \to \infty$ for fixed $t>0$ is impossible and
the conditions of the CLT in Theorem \ref{thm:jumps_vanish} are not satisfied either when $\beta>1$.
If  $\beta<1$ then $\langle Z_n \rangle_t \to 0$ holds and so do all the conditions of this CLT.
This leaves only the case $\beta=1$ and $\alpha=\delta>0$: neither condition \eqref{eq:condM} nor
the CLT in Theorem \ref{thm:jumps_vanish} hold.
Summarizing: under $\alpha>0$ and $\delta>0$, the Lindeberg condition \eqref{eq:condL}, the McLeish condition \eqref{eq:condM}
and those for the CLT in Theorem \ref{thm:jumps_vanish} all hold when $\beta<1$ but do not when $\beta\ge1$.
For this example the three are equivalent.
The construction of a sequence of martingales $B_n$ converging to a Brownian motion $B$ with the desired properties
in order to apply Theorem \ref{thm:jumps_vanish}, has been sidestepped here.
Enlarge the space (if necessary) in order to use an iid Bernoulli sequence $\{\psi_i\}$ with the right features.
Define the martingale $\disp B_n(t):=\frac{1}{\sn}\sum_{k=1}^\nt \psi_k \sign(\xi_k)$, with
$\{\psi_k\}$ iid with $\disp P(\psi_k=\pm 1) = \frac{1}{2}$, independent of the whole sequence $\{X_j\}$.
That $B_n \stackrel{\mathcal{C}}{\rightsquigarrow} B$ holds for some Brownian motion $B$ results from another application of
Theorem \ref{thm:jumps_vanish} through the observations : $\langle M_n,B_n  \rangle = [M_n,B_n]_t = 0$,
$\disp \crochet{ B_n}_t=\frac{\nt}{n}\to t$ and the use of a further martingale $\disp \frac{1}{\sn}\sum_{k=1}^\nt \psi_k$.
Furthermore, while uniform convergence of $A_n$ fails when $\delta=0$, it is restored by applying Theorem \ref{thm:jumps_vanish} to
rescaled martingale $\disp \tilde M_n(t) = \{\log{m(n)}\}^{-1/2} \sum_{k=1}^{\lfloor m(nt) \rfloor} \xi_k$ with the
normalization $a_n^2=\log{m(n)}$ where $m(t) = te^t$ is used instead of $m(t)=t$,
maintaining the requirement $\tilde M_n(0)=0$ and finding
$(\tilde M_n,\langle \tilde M_n \rangle) \stackrel{\mathcal{C}}{\rightsquigarrow} (\tilde M,\langle \tilde M \rangle)$
with $\langle \tilde M \rangle_t=t$ and hence $\tilde M$ a Brownian motion.
Note also that, for this rescaled martingale, when $\delta=0$ and $\beta>1$,
$[\tilde M_n]_t \stackrel{a.s.}{\rightsquigarrow} 0\neq[\tilde M]_t$ and $[\tilde M_n] \stackrel{\mathcal{C}}{\rightsquigarrow} 0\neq[\tilde M]$,
illustrating the discontinuity of the map $\tilde M\mapsto[\tilde M]_t$ for any fixed $t>0$, in general.
Finally, when $\alpha=\delta>0$ and $\beta=1$, the Lindeberg condition \eqref{eq:condL},
the McLeish condition \eqref{eq:condM} and the conditions of the CLT in Theorem \ref{thm:jumps_vanish} all fail;
however, Corollary \ref{cor:main_terr1b} is applicable here since $A_n(t)=E\{[M_n]_t\}$.
\end{exmp}


\subsection{Applications}\label{ssec:linnik_gamma}

The following example, inspired from Mathematical Finance, shows that one can indeed get discontinuous limits.

\begin{exmp}\label{exmp:linnik}
Let $\{\xi_{n,k}:k\ge1,n\ge1\}$ be independent ${\rm Gamma}(1/n,1)$, independent of $\{Z_k:k\ge1\}$ which are iid ${\rm N}(0,1)$.
Set $\disp M_n(t)=\sum_{k=1}^\nt \xi_{n,k}^{1/2}Z_k$ and define filtration $\cF_{n,k} = \sigma\{Z_1,\ldots, Z_k\}\cup\cF_{n,0}$,
with initial $\sigma$-algebra $\cF_{n,0} = \sigma\{\xi_{n,k},k\ge1\}$.
Then $M_n$ is an $\cF_{n,[n\cdot]}$-martingale with quadratic variation $[M_n]_t = \sum_{k=1}^\nt \xi_{n,k}Z_k^2$
and predictable compensator $\disp A_n(t) = \sum_{k=1}^\nt \xi_{n,k}$. First there holds
$A_n \stackrel{f.d.d.}{\rightsquigarrow} A$, where $A(t)\sim {\rm Gamma}(t,1)$.
We also have $M_n \stackrel{f.d.d.}{\rightsquigarrow} M$,
where $M \stackrel{Law}{=} W\circ A$, and $W$ is a Brownian motion independent of $A$. In fact, $M(t)$ has a Linnik distribution
\citep{Jacques/Remillard/Theodorescu:1999} with characteristic function $E\left\{e^{i\lambda M(t)}\right\}=(1+\lambda^2/2)^{-t}$.
Since $\esp{A_n(t)-A_n(s)}=\esp{|M_n(t)-M_n(s)|^2}=\disp\frac{\nt-\ns}{n}$, both
$A_n\stackrel{\cJ_1}{\rightsquigarrow} A$ and $M_n\stackrel{\cJ_1}{\rightsquigarrow} M$ hold,
using \cite[Theorem]{Genest/Ghoudi/Remillard:1996} on each time interval $[0,T]$.
The rest of Hypothesis \ref{hyp:An2} holds --- for instance, Hypothesis \ref{hyp:An2}d comes from
$\disp \crochet{ W_n }_t = A_n\circ \tau_n(t) = V_{n,V_n^{-1}(t)} \in\left[t,t+\xi_{n,V_n^{-1}(t)}\right]$.
Since $A$ is a strictly increasing L\'evy process, $\tau$ is continuous everywhere and therefore
Remark \ref{rem:strict} implies that $\tau$ itself is not a L\'evy process. Nevertheless, all the requirements for the use of
Corollary \ref{cor:main_terr1b} are met, yielding $(M_n,A_n) \stackrel{\cJ_1}{\rightsquigarrow} (M,A)$.
\end{exmp}


 \begin{rem}\label{rem:levy}
The previous limit is a special case of the so-called Variance Gamma process \citep{Madan/Carr/Chang:1998}, which appears as models for log-return prices. See also \cite{Ane/Geman:2000}. Note that for any subordinator $L$ (i.e., increasing L\'evy process) with finite mean, one can replace $\xi_{n,k}$ by $L(k/n)-L((k-1)/n)$, $k\ge 1$. In the limit one obtains $M = W\circ L$, where $W$ is a Brownian motion independent of $L$. Another popular subordinator is the Inverse Gaussian process, leading to the so-called Normal Inverse Gaussian process \citep{Remillard:2013}.
\end{rem}

 \begin{exmp}\label{exmp:anegeman}
For modeling log-returns,  \citet{Ane/Geman:2000} introduced the model $W\circ A+\mu A$, i.e., a time-changed BM with drift,
for the (adjusted) financial return for investments, taking into account a time scale $A$ which is meant to reflect business cycles and other features known collectively in economics as business time. A detailed analysis of the asymptotic properties of realized volatility error in investment returns, for the special case where $A$ has continuous trajectories, was first obtained in \cite{Barndorff-Nielsen/Shephard:2002} and extended in \citet{Remillard/Vaillancourt:2024a} --- notably, proofs of both the consistency and the invariance (up the order of third moments) for their now popular estimator of the realized volatility error in investment returns, when data is collected at regular intervals.
When $A$ is discontinuous, a more realistic approach to financial fluctuations, take another look at Example \ref{exmp:linnik}.
Martingales $\disp M_n(t)=\sum_{k=1}^\nt \xi_{n,k}^{1/2}Z_k$ with respective compensators $\disp A_n(t) = \sum_{k=1}^\nt \xi_{n,k}$ are shown to verify $(M_n,A_n) \stackrel{\cJ_1}{\rightsquigarrow} (W\circ A,A)$, where $W$ is a Brownian motion independent of $A$. Since the convergence is joint, there ensues at once $M_n+\mu A_n \stackrel{\cJ_1}{\rightsquigarrow} W\circ A+\mu A$!
In the light of Remark \ref{rem:levy}, this construction yields such financial return models for all L\'evy subordinators for Brownian motion.
  \end{exmp}


\section{Non-stationary arrays with discontinuous limits}\label{sec:MPU}

Given is an array $\{Y_{n,j}\}_{j,n\ge 1}$ of random variables with mean $0$
and finite variances, not necessarily all the same. No restriction is made on their dependency structure.
Put $\cG_{n,k} := \sigma\{Y_{n,j}; j\le k\}$ and successively set $Z_{n,j}:=E\left(Y_{n,j}|\cG_{n,j-1}\right)$, $X_{n,j}:=Y_{n,j}-Z_{n,j}$,
$\disp N_{n}(t) := \sum_{j=1}^{\nt} Y_{n,j}$ and $\disp O_{n}(t) := \sum_{j=1}^{\nt} Z_{n,j}$. Note that $\cG_{n,k}=\sigma\{X_{n,j}; j\le k\}$ as well.
The square integrable martingale of interest here is $\disp M_{n}(t)=N_{n}(t)-O_{n}(t)= \sum_{j=1}^{\nt} X_{n,j}$,
with $A_n(t) = \langle M_n\rangle(t)=\disp \sum_{k=1}^\nt E\left(X_{n,k}^2|\cG_{n,k-1}\right)$.
Empty sums are null, therefore in particular set $N_n(0)=O_n(0)=M_n(0)=A_n(0)=0$.
\citet{Merlevede/Peligrad/Utev:2019} contains several CLTs of the form
$N_n \stackrel{\mathcal{C}}{\rightsquigarrow} N_\infty$
when the limit $N_\infty$ is the stochastic integral of a continuous function with respect to a Brownian motion.
Corollary \ref{cor:main_terr1b} allows for the extension of several of their results to mixtures of stochastic
processes including some with discontinuities. Typically, Corollary \ref{cor:main_terr1b} together with
Lemma \ref{lem:BasicInequality} yield $N_n \stackrel{\cJ_1}{\rightsquigarrow} N_\infty = W\circ A+O_\infty$ with
$W$ is a Brownian motion independent of $A$, provided all conditions are verified of course, notably
$O_n \stackrel{\mathcal{C}}{\rightsquigarrow} O_\infty$ with a continuous limit.
The following example gives a sample calculation of this type of asymptotic transfer.

\begin{exmp}\label{exmp:burkholder}
Consider an arbitrary array $\{Q_{n,j}\}_{j,n\ge 1}$ of
bounded predictable random variables in the sense that each $Q_{n,j}$ is $\cG_{n,j-1}$-measurable.
Random variables $Y_{n,j}$ still have mean $0$ and finite second moment.
Set $\bar N_{n}(t) := \sum_{j=1}^{\nt} Q_{n,j}Y_{n,j}$, a weighted version of $N_n$,
$\disp \bar O_{n}(t) := \sum_{j=1}^{\nt} Q_{n,j}Z_{n,j}$ and
$\disp \bar M_{n}(t) := \sum_{j=1}^{\nt} Q_{n,j}X_{n,j}$.
The martingale transform $\bar M_{n}$ of $\{X_{n,j}\}_{j,n\ge 1}$ by $\{Q_{n,j}\}_{j,n\ge 1}$
\citep{Burkholder:1966} is another martingale, with compensator
$\disp \langle \bar M_n \rangle _t = \sum_{j=1}^{\nt}Q_{n,j}^2E\left(X_{n,j}^2|\cG_{n,j-1}\right)$.
Set $\disp\bar Q_n(t):=\sum_{j=1}^\infty Q_{n,j}\I_{\{j/n\le t<(j+1)/n\}}$, $\bar W_n:=\bar M_n\circ\bar \tau_n$
and $\langle \bar W_n\rangle:=\langle \bar M_n\rangle\circ\bar\tau_n$.
Take notice that $\langle \bar M_n \rangle _t = \int_0^t \bar Q_n^2(s) dA_n(s)$ and
$\bar O_{n}(t) = \int_0^t \bar Q_n(s) dO_n(s)$, both of which are Lebesgue-Stieltjes integrals.

The following result provides a prototypical CLT obtained by such an asymptotic transfer.
In essence, it is a CLT for semimartingales. It is proved in Appendix \ref{pf:thm:mart_trans}.

\begin{thm}\label{thm:mart_trans}
Assume that $(M_n,A_n,\tau_n,W_n)$ and their limits satisfy all the hypotheses in the first sentence of Corollary \ref{cor:main_terr1b}.
Assume also all of the following:
\begin{itemize}
\item[a)] $O_n \stackrel{\mathcal{C}}{\rightsquigarrow} O_\infty$ holds for some process $O_\infty$ with bounded variations on
compact time sets;

\item[b)] $\disp 0<\inf_{j,n\ge1}|Q_{n,j}|\le \sup_{j,n\ge1}|Q_{n,j}|<\infty$;

\item[c)] $\bar Q_n \stackrel{\mathcal{C}}{\rightsquigarrow} \bar Q_\infty$ holds for some continuous process $\bar Q_\infty$;

\item[d)] $\bar A_\infty(t):=\int_0^t \bar Q_\infty^2(u)dA(u)$ is an inhomogeneous L\'evy process.
\end{itemize}
Then $\bar N_n \stackrel{\cJ_1}{\rightsquigarrow} \bar N_\infty=\bar W\circ \bar A_\infty+\bar O_\infty$ holds with
$\bar O_\infty(t) = \int_0^t \bar Q_\infty(s) dO_\infty(s)$ and $\bar W$ a Brownian motion independent of $\bar A_\infty$.
\end{thm}

Thus, under some technical restrictions, the two simpler CLTs a) $O_n \stackrel{\mathcal{C}}{\rightsquigarrow} O_\infty$ and
c) $\bar Q_n \stackrel{\mathcal{C}}{\rightsquigarrow} \bar Q_\infty$ can be combined to yield the sought after CLT
$\bar N_n \stackrel{\cJ_1}{\rightsquigarrow} \bar N_\infty$ with a discontinuous limit. Here is a sample calculation.
 \end{exmp}

\begin{exmp}\label{exmp:jacquesetal1}
Pursuant to the notation introduced in Example \ref{exmp:linnik}, recall that martingales
$\disp M_n(t)=\sum_{k=1}^\nt \xi_{n,k}^{1/2}Z_k$ with respective compensators $\disp A_n(t) = \sum_{k=1}^\nt \xi_{n,k}$
verify $(M_n,A_n) \stackrel{\cJ_1}{\rightsquigarrow} (M,A)$, essentially through Corollary \ref{cor:main_terr1b}.
Let $Q$ denote a continuous process on $[0,1]$, valued in some compact subset of $(-\infty,0)\cup(0,\infty)$
in order to avoid trivialities related to null or near null weights.
Assume $Q$ is independent of $\{\xi_{n,k}:k\ge1,n\ge1\}\cup\{Z_k:k\ge1\}$.
Set $\disp Q_{n,k}=Q\left(\frac{k-1}{n}\right)$. Think of $Q$ as a random weight representing some fluctuating linear filter
superimposed on the underlying motion $M_n$. Set $\disp \bar M_{n}(t) := \sum_{k=1}^{\nt} Q_{n,k}\xi_{n,k}^{1/2}Z_k$
(Example \ref{exmp:burkholder}). Using $O_n=O_\infty \equiv 0$, Theorem \ref{thm:mart_trans} yields
$\bar M_n \stackrel{\cJ_1}{\rightsquigarrow} \bar W\circ \bar A_\infty$ with $\bar Q_\infty=Q$,
since $\bar A_\infty(t):=\int_0^t Q^2(u)dA(u)$ is an inhomogeneous L\'evy process,
by \citet[Theorem II.4.15 and Proposition IX.5.3]{Jacod/Shiryaev:2003}.
Furthermore, for any sequence $\{O_n\}$ and limit $O_\infty$ satisfying condition a), there also holds
$\bar N_n \stackrel{\cJ_1}{\rightsquigarrow} \bar W\circ \bar A_\infty+\bar O_\infty$.
 \end{exmp}

\begin{exmp}\label{exmp:jacquesetal2}
Alternatively to Example \ref{exmp:jacquesetal1}, Example \ref{exmp:linnik} can also be extended in the following manner.
Set $\disp \phi_{n,i}=n^{-1/2}\phi_{i}$ where $\{\phi_{i}:i\ge0\}$ are centered iid random variables with finite second moments
and independent of $\{\xi_{n,k}\}\cup\{Z_k\}$.
Let $Q$ denote a continuous function on $[0,1]$, valued in some compact subset of $(-\infty,0)\cup(0,\infty)$.
Put $\disp Q_{n,k}=Q\left(\sum_{i=0}^{k-1}\phi_{n,i}\right)$ and set
$\disp \bar M_{n}(t) := \sum_{k=1}^{\nt} Q_{n,k}\xi_{n,k}^{1/2}Z_k$ (again using Example \ref{exmp:burkholder}).
Hypotheses b) and c) of Theorem \ref{thm:mart_trans} hold with $\bar Q_\infty=Q(\sigma B)$, where $\sigma^2:=E\{\phi_i^2\}$
and $B$ is a Brownian motion independent of $(W,A)$. There ensues
$\bar M_n \stackrel{\cJ_1}{\rightsquigarrow} \bar W\circ \bar A_\infty$
since $\bar A_\infty(t):=\int_0^t Q^2(\sigma B_u)dA(u)$ is also an inhomogeneous L\'evy process, actually a special case
of the one encountered in Example \ref{exmp:jacquesetal1}.
Here again, there also holds $\bar N_n \stackrel{\cJ_1}{\rightsquigarrow} \bar W\circ \bar A_\infty+\bar O_\infty$
for any sequence $\{O_n\}$ and limit $O_\infty$ satisfying condition a).
 \end{exmp}


\appendix

\section{Real-valued L\'evy  processes}\label{app:auxresultsLevy}

An inhomogeneous L\'evy process $A$ is an $\mathbb{F}$-adapted $D$-valued process with independent increments (with respect to
$\mathbb{F}$), without fixed times of discontinuity and such that $A_0=0\in\dR$. The inhomogeneity (in time) means that stationarity of the increments, usually required of L\'evy processes, is lifted. This choice reflects the existence of weak limits exhibiting independent increments without homogeneity in time or space, in some applications. All processes are built directly on $D$, ensuring continuity in probability of the trajectories, another usual requirement.

The characteristic function of an inhomogeneous L\'evy process $A$ is given by
$$
E\left[\left. e^{i\theta (A_t-A_s)}\right| \cF_s \right] = e^{\Psi_t(\theta)-\Psi_s(\theta)}, \quad 0\le s < t,
$$
where $A$ is $\mathbb{F}$-adapted with $\dF = (\cF_t)_{t\ge 0}$,
$$
\Psi_t(\theta) = i \theta B_t  -\frac{\theta^2 C_t}{2}
+ \left\{\prod_{0<r\le t} e^{-ir\Delta B_r}\right\}\int_0^t \int \left(e^{i\theta x}-1-i\theta x\I_{\{|x|\le 1\}}\right)\nu(du,dx),
$$
$\Delta B_r = B_r-B_{r-}$, and $\nu$ is a L\'evy measure.
For details see \citet[Theorem II.4.15, Theorem II.5.2 and Corollaries II.5.11,12,13]{Jacod/Shiryaev:2003}.
The characteristics of $A$ are defined to be the deterministic functions $(B,C,\nu)$.
The choice of $h(x)=x\I_{\{|x|\le 1\}}$ is classical but in \citet[Remark II.2.7]{Jacod/Shiryaev:2003},
it is rejected in favor of a continuous function $h$, a substitution which affects the value of $B$ exclusively
--- see \citet[Corollary II.5.13]{Jacod/Shiryaev:2003}.
One can instead make use of $h(x)=x\I_{\{|x|\le 1\}}+\sign(x)\I_{\{|x|>1\}}$, which is continuous.
A L\'evy process is called homogeneous or stationary when $B_t=tB_1$, $C_t=tC_1$ and $\nu([0,t]\times A)=t\nu([0,1]\times A)$,
so that $\Psi_t(\theta) = t\Psi_1(\theta)$ ensues.
Important stationary L\'evy processes include Brownian motion ($\nu\equiv0$) and, for each $\alpha\in(0,2)$,
(non-Gaussian) $\alpha$-stable L\'evy motions \citep[Section 4.5]{Whitt:2002}, characterized by both $C_t\equiv0$ and
$\disp \nu(du,dx)=dudx\cdot|x|^{-1-\alpha}\left[c_\alpha^{+}\I_{\{x>0\}}+c_\alpha^{-}\I_{\{x<0\}}\right]$, when $\alpha\in(0,1)\cup(1,2)$,
with an alternative form for $\nu$ when $\alpha=1$, the Cauchy process, which will not concern us in this paper.
Recall that $\alpha$-stable L\'evy motions are either without negative jumps ($1\le\alpha<2$) or non-decreasing ($0<\alpha<1$).

The following key result is taken from \citet{Remillard/Vaillancourt:2025}.

\begin{lem}\label{lem:main}
Suppose $A$ is an $\mathbb{F}$-adapted inhomogeneous L\'evy process with finite variation on any finite interval $[0,T]$ and no Brownian component, i.e., $C_t = 0$ for all $t\ge 0$. Then $A$ is independent of any $\mathbb{F}$-Brownian motion $W$.
Further, if $A$ is non-negative and non-decreasing everywhere, then the random rescaling property
\begin{equation}\label{eq:random_rescaling}
W\circ A (t) - W\circ A (s) \stackrel{Law}{=} (W_t-W_s)\sqrt{\frac{A_t-A_s}{t-s}}
\end{equation}
holds for every $0\le s\le t<\infty$ and both sides define equidistributed inhomogeneous L\'evy processes.
\end{lem}


\section{Topologies and weak convergence on $D$}\label{app:M1Topology}

Equality in law is denoted by $\stackrel{Law}{=}$,
convergence in law by $\stackrel{Law}{\rightsquigarrow}$,
in probability by $\stackrel{Pr}{\rightsquigarrow}$ and
almost sure convergence by $\stackrel{a.s.}{\rightsquigarrow}$.
Sequences of random processes are sometimes written $(X_n)_{n\ge 1}$ or $\{X_n\}_{n\ge 1}$ for clarity,
but usually just as $X_n$ for brevity; similarly for all other types of sequences.

For either $\dH=\dJ$ or $\dH=\dM$, a (deterministic) sequence $x_n\in D$ is said to
$\mathcal{H}_1$-converge to $x\in D$ if and only if
\begin{enumerate}
\item[(i)]
$\lim_{n\to\infty}x_n(t)=x(t)$ holds pointwise for every $t$ in an everywhere dense subset of $[0,\infty)$ that includes the origin,

\item[(ii)]
 for every $T>0$,
\begin{equation}\label{eq:modulus1}
\lim_{\delta\to 0}\limsup_{n\to\infty}\omega_\dH(x_n,\delta,T) =0.
\end{equation}
\end{enumerate}
A subset $K\subset D$ is $\mathcal{H}_1$-compact if and only if, for every $T>0$, there holds
$\sup_{x\in K}\sup_{0\le t\le T}|x(t)|<\infty$ and
$\lim_{\delta\to 0}\sup_{x\in K}\omega_\dH(x,\delta,T) =0$.

A sequence of $D$-valued processes $X_n$ are $\mathcal{H}_1$-tight if and only if
\begin{enumerate}
\item[(i)]
for every $t$ in an everywhere dense subset of $[0,\infty)$ that includes the origin, the marginal distributions of $X_n(t)$ are tight,

\item[(ii)]
 for any $\epsilon>0$, and $T>0$,
\begin{equation}\label{eq:modulus2}
\lim_{\delta\to 0}\limsup_{n\to\infty}P\{\omega_\dH(X_n,\delta,T) > \epsilon\}=0.
\end{equation}
\end{enumerate}
Recall that $\mathcal{C}$-tightness is defined by setting $\dH=\dC$
with $\dC(x_1,x_2,x_3):=|x_3-x_2|$ in (\ref{eq:modulus2}).
When (\ref{eq:modulus2}) holds for any $T$, it remains so on the boundary $T=\delta$,
for each $\dH=\dC$, $\dH=\dJ$ or $\dH=\dM$.

Complete coverage of the $\cJ_1$-topology can be found in \citet{Ethier/Kurtz:1986} and \citet{Jacod/Shiryaev:2003},
where $\dB$ is actually allowed be a Polish space; while \citet{Whitt:2002} is our main reference on the $\cM_1$-topology.
Neither topology turns $D$ into a topological vector space, that is to say, neither is compatible with the linear structure inherited
on $D$ from the Banach space $\dB$, which would have made the sum a continuous operator. Consequently, when considering vectors of processes valued in a product $\prod_{i=1}^d \dB_i$ of Banach spaces, neither topology offers the equivalence
between $\mathcal{H}_1$-convergence of sequences of probability measures on $D$ and coordinatewise $\mathcal{H}_1$-convergence,
an equivalence which holds for $\mathcal{C}$-convergence.
Similarly $\mathcal{H}_1$-tightness cannot be treated coordinatewise; it does however ensue from the
$\mathcal{H}_1$-tightness of each coordinate plus that of each pairwise sum of coordinates --- $d^2$ verifications are thus required.
This fact precludes immediate extensions of some results from the real line to higher dimensional spaces, but there are exceptions,
as will be seen in several of the proofs in Appendix \ref{app:main_results}. One such exception is the following result, taken from \citet{Remillard/Vaillancourt:2025}.

\begin{lem}\label{lem:BasicInequality}
For any pair of functions $X, Y\in D$ there holds, for every choice of $\delta>0$ and $T>0$,
$\omega_\dH(X+Y,\delta,T)\le \omega_\dH(X,\delta,T)+\omega_\dC(Y,\delta,T)$.
For sequences of $D$-valued processes, if $X_n$ is $\mathcal{H}_1$-tight and
$Y_n$ is $\mathcal{C}$-tight, then $X_n+Y_n$ and $(X_n,Y_n)$ are also $\mathcal{H}_1$-tight.
Similarly if $X_n\stackrel{\mathcal{H}_1}{\rightsquigarrow} X$, $Y_n\stackrel{\mathcal{C}}{\rightsquigarrow} Y$, and
$(X_n,Y_n)\stackrel{f.d.d.}{\rightsquigarrow} (X,Y)$ then $X_n+Y_n\stackrel{\mathcal{H}_1}{\rightsquigarrow} X+Y$
and $(X_n,Y_n)\stackrel{\mathcal{H}_1}{\rightsquigarrow}(X,Y)$.
\end{lem}

Deducing the $\mathcal{H}_1$-tightness of $D$-valued processes from that of dominating non-decreasing processes,
is achieved using Lenglart's celebrated inequality, the proof of which can be found in
\citet[Lemma I.3.30]{Jacod/Shiryaev:2003}.
For any $f \in D$ and $T>0$, set $\Delta f(t) := f(t)-f(t-)$ and $J_T(f) := \sup_{t\in [0,T]}|\Delta f(t)|<\infty$.

\begin{lem}[Lenglart's inequality]\label{lem:lenglart}
Let $X$ be an $\dF$-adapted $D$-valued process. Suppose that $Y$ is optional, non-decreasing,
and that, for any bounded stopping time $\tau$, $E|X(\tau)|\leq \esp{Y(\tau)}$.
Then for any stopping time $\tau$ and all $\varepsilon,\eta >0$,
\begin{itemize}
\item[{a)}]
if $Y$ is predictable,
\begin{equation}\label{eng1}
P(\sup_{s\leq \tau}|X(s)|\geq\varepsilon) \leq \frac{\eta}{\varepsilon} +
P(Y(\tau)\geq\eta).
\end{equation}

\item[b)]
if $Y$ is adapted,
\begin{equation}\label{eng2}
P(\sup_{s\leq \tau}|X(s)|\geq\varepsilon) \leq \frac{1}{\varepsilon}
\left[\eta+ E\left\{J_\tau(Y) \right\}\right] +
P(Y(\tau)\geq\eta).
\end{equation}
\end{itemize}
\end{lem}

Proving the stronger $\cJ_1$-tightness generally involves the following lemma.

\begin{lem}[Aldous's criterion]\label{lem:aldous}
Let $\{X_n\}_{n\ge 1}$ be a sequence of $D$-valued processes. Suppose that for any sequence of bounded discrete stopping times
$\seq{\tau_n}$ and for any sequence $\seq{\delta_n}$ in $[0,1]$ converging to $0$, the following condition holds, for every $T>0$:
(A) $X_n((\tau_n+\delta_n)\wedge T)- X_n(\tau_n)\stackrel{Law}{\rightsquigarrow} 0$. Then, $\seq{X_n}$ is $\cJ_1$-tight,
if either of the two following conditions holds:
\begin{enumerate}
\item $\{X_n(0)\}_{n\ge 1}$ and $(J_T(X_n))_{n\ge 1}$ are tight;
\item $\{X_n(t)\}_{n\ge 1}$ is tight for any $t\in [0,T]$.
 \end{enumerate}

\end{lem}

\proof See \citet{Aldous:1978}, \citet[Theorem VI.4.5]{Jacod/Shiryaev:2003}
or for several variants see \citet[Theorem 3.8.6]{Ethier/Kurtz:1986}.
\qed

\begin{rem}\label{rem:aldous}
Note that condition (1) or condition (2) are necessary for $\cJ_1$-tightness, but not condition (A).
Also note that condition (A) holds when $X_n$ is $\mathcal{C}$-tight, since
$|X_n((\tau_n+\delta_n)\wedge T)- X_n(\tau_n)|\leq \omega_\dC \left(X_n,\delta_n,T\right)$.
 \end{rem}


\section{Non-decreasing functions and processes}\label{app:continuousinverse}

In this section we focus exclusively on the real-valued functions and processes so $\dB=\dR$ throughout.

For any c\`adl\`ag non-decreasing non-negative function $A\in D$, denote its (compositional) inverse by $\tau(s) = \inf\{t\ge0; A(t)>s\}$.
Let $D_{\uparrow}\subset D$ be the subspace of those non-decreasing non-negative $A$ such that $A(t) \to\infty$ as $t\to\infty$
and let $D_{\uparrow}^0\subset D_{\uparrow}$ be the further subspace where $\tau(0)=0$ as well.
The following result, on the $\cM_1$-continuity of the inverse map (which in general fails to be $\cJ_1$-continuous),
is from \citet{Remillard/Vaillancourt:2025}.

\begin{prop}\label{prop:inversecontinuity}
The inverse map $A\mapsto\tau$ is a well defined bijective mapping from $D_{\uparrow}$ into itself, such that there holds
$A\circ\tau\circ A=A$ provided $A\in D_{\uparrow}$ is either continuous everywhere or strictly increasing everywhere.
Both $A\mapsto\tau$ and the reverse map $\tau\mapsto A$ are $\cM_1$-continuous
when restricted to $D_{\uparrow}^0$; this is not necessarily the case without this additional restriction, not even on $D_{\uparrow}$.
\end{prop}

The following result is from \citet[Remark 3.1]{Remillard/Vaillancourt:2025}.

\begin{prop}\label{prop:Cincreasing}
Let $D$-valued non-decreasing processes $A_n$ and some continuous process $A$
be such that $A_n \stackrel{f.d.d.}{\rightsquigarrow} A$. Then $A_n$ is $\mathcal{C}$-tight
and $A_n \stackrel{\mathcal{C}}{\rightsquigarrow} A$.
\end{prop}


\section{Tightness}\label{app:TightnessTheMainLemma}

In this section the martingales take their values in $\dB=\dR$.

Let $M$ be a $D$-valued square integrable $\mathbb{F}$-martingale started at $M(0)=0$ and
with quadratic variation $[M]$, predictable compensator $\langle M \rangle $, continuous part $M^c$,
purely discontinuous part $M^d$ and largest jump $J_T(M)$ up to time $T>0$.
The existence of $M^c$ and $M^d$ as orthogonal square integrable $\mathbb{F}$-martingales
such that $M=M^c+M^d$ was established during the early development of martingale theory
and a proof can be found in \citet[Theorem I.4.18]{Jacod/Shiryaev:2003}.

\begin{prop}\label{prop:decomp}
Both the following hold almost surely:
\begin{eqnarray*}
[M]_T &=& \langle M^c \rangle _T + \sum_{0<s\le T}\{\Delta M(s)\}^2;\\
\langle M \rangle _T &=& \langle M^c \rangle _T + \sum_{0<s\le T}E\bigl(\{\Delta M(s)\}^2|{\cF}_{s-}\bigr).
\end{eqnarray*}
\end{prop}
\proof
The first statement is actually valid for any semimartingale $M$
with $M^c$ denoting its continuous martingale part \citep[Theorem I.4.52]{Jacod/Shiryaev:2003}.
The second one is equivalent to stating that
$$
N_T = \sum_{0\le s \le T} \Delta M^2(s) - \sum_{0\le s \le T} E\left\{ \Delta M^2(s)|\cF_{s-}\right\}
$$
is a martingale. This is indeed true. For if $t<T$, then
\begin{eqnarray*}
E(N_T|\cF_t) &=& N_t + E\left[ \left.
\sum_{t< s\le T} \Delta M^2(s)  -  E\left\{  \Delta M^2(s) | \cF_{s-} \right\} \right| \cF_t\right] \\
& = & N_t + \sum_{t< s\le T}  \left[ E\left\{ \Delta M^2(s)|\cF_t\right\} - E\left\{ \Delta M^2(s)|\cF_t \right\}\right] =N_t,
\end{eqnarray*}
because $E \left[ \E\left\{ \Delta M^2(s) |\cF_{s-}\right\}|\cF_t\right] =  E\left\{ \Delta M^2(s)|\cF_t \right\}$ for any $t<s$.
\qed

Note that $\Delta \langle M \rangle _s = \esp{\Delta [M]_s|{\cF}_{s-}}$
and $\esp{J_s(\langle M \rangle)}\le \esp{J_s([M])}$, for further reference.

Let $M_n$ be a sequence of $D$-valued square integrable $\mathbb{F}$-martingales started at $M_n(0)=0$
with quadratic variation $[M_n]$ and predictable compensator $\langle M_n \rangle $.


Now comes the main result about $\cJ_1$-tightness, in the real-valued case $\dB=\dR$.

\begin{thm}\label{thm:M1tightness}
Let $M_n$ be a sequence of $D$-valued square integrable $\mathbb{F}$-martingales
with predictable quadratic variation $A_n = \langle M_n \rangle $, such that $(M_n(0))_{n\ge 1}$ is tight.

\begin{itemize}
\item[a)]
If $A_n$ verifies condition (A) in Lemma \ref{lem:aldous}, then so do both $M_n$ and $[M_n]$.
If in addition $J_t(M_n)$ is tight for every $t>0$, then both $M_n$ and $[M_n]$ are $\cJ_1$-tight.
Finally, if $J_t(M_n)\stackrel{Law}{\longrightarrow}0$, then both $M_n$ and $[M_n]$ are $\mathcal{C}$-tight.

\item[b)]
Assume $\limsup_{n\to\infty} \esp{J_t^2(M_n)}<\infty$ for any $t>0$.
If $[M_n]$ verifies condition (A) in Lemma \ref{lem:aldous}, then so do both $M_n$ and $A_n$; all three are $\cJ_1$-tight.
If in addition $J_t(M_n)\stackrel{Law}{\longrightarrow}0$, then all three are $\mathcal{C}$-tight.

\item[c)]
If $A_n$ is $\mathcal{C}$-tight, then both $M_n$ and $[M_n]$ are $\cJ_1$-tight.

\item[d)]
Assume $\sup_nE\{J_t(M_n)\}<\infty$ for every $t>0$. If $M_n$ is $\cJ_1$-tight then so is $[M_n]$.
\end{itemize}
\end{thm}

\begin{rem}\label{rem:classic}
$\mathcal{C}$-tightness implies condition (A) in Lemma \ref{lem:aldous} (Remark \ref{rem:aldous}).
Even the $\mathcal{C}$-tightness of $[M_n]$ does not necessarily imply the tightness of $A_n(t)$ for a single $t>0$ without some
control over the growth of the latter  --- choose scaling $a_n^2=n^2$ and any pair $(\alpha,\beta)$ satisfying $\alpha-1>\beta>1$
in Example \ref{exmp:lindeberg}, where $[M_n]\stackrel{\mathcal{C}}{\rightsquigarrow}0$ while $\lim_nA_n(t)=\infty$ for every $t>0$.
When the increments of $M_n$ are uniformly bounded, the tightness of $[M_n]_t$ for all $t>0$ and that of $A_n(t)$ for all $t>0$
are equivalent and each imply that of $M_n(t)$ for all $t>0$ as well \citep[Proposition VI.6.13]{Jacod/Shiryaev:2003}.
For additional insight into the role Lemma \ref{lem:aldous}
plays in the proof of $\cJ_1$-tightness results, see \citet[Remark VI.4.7]{Jacod/Shiryaev:2003}.
 \end{rem}

\proof a) Set $X_n(s)=\left\{M_n (s+\tau_n)-M_n(\tau_n)\right\}^2$ and
$Y_n(s)= \langle  M_n  \rangle _{s+\tau_n}-\langle  M_n  \rangle _{\tau_n}$,
where $\tau_n$ is a sequence of stopping times uniformly bounded by $T$ for any $n$,
and $\esp{X_n(\tau)} = \esp{Y_n(\tau)}$, for any bounded stopping time $\tau$.
Using Equation \eqref{eng1} from Lemma \ref{lem:lenglart}, we have, for any $\varepsilon>0$,
$\eta>0$ and any numerical sequence $\delta_n \in (0,1) \to 0$,
\begin{equation*}
P\left\{ |M_n(\tau_n+\delta_n)-M_n(\tau_n)|\geq\varepsilon \right\}  \leq
P\left\{Y_n(\delta_n)>\eta\right\}  + \frac{\eta}{\varepsilon^2}.
\end{equation*}
Choosing $\eta = \varepsilon^3$,
$
\limsup_{n\to\infty}
P\left\{ |M_n(\tau_n+\delta_n)-M_n(\tau_n)|\geq\varepsilon \right\} \leq \varepsilon
$
follows, since $Y_n(\delta_n)$ tends to $0$ as $n\to\infty$, showing that condition (A) of Lemma \ref{lem:aldous} is met by $M_n$ as well.
When condition (1) is assumed, $M_n$ is $\cJ_1$-tight.
Finally, if $J_t(M_n)\stackrel{Law}{\longrightarrow}0$, then $\mathcal{C}$-tightness and $\cJ_1$-tightness are equivalent for $M_n$.
The proof is similar for $[M_n]$, using $J_T([M_n])=J_T^2(M_n)$.

b) Set $X_n(s) =  \langle  M_n  \rangle _{s+\tau_n}-\langle  M_n  \rangle _{\tau_n}$
and $Y_n(s) = [M_n]_{s+\tau_n}-[M_n]_{\tau_n}$ instead, so there holds
$E \{X_n(\delta_n)\} = E \{Y_n(\delta_n)\}$ and
$$
P(\langle  M_n  \rangle _{s+\tau_n}-\langle  M_n  \rangle _{\tau_n} \geq\varepsilon) \leq \frac{1}{\varepsilon}
\left[\eta+ E\left\{J_{\delta_n}(Y_n) \right\}\right] + P(Y_n(\delta_n))\geq\eta),
$$
using \eqref{eng2} of Lemma \ref{lem:lenglart}, with $J_{\delta_n}(Y_n)\le Y_n(\delta_n)$ since $[M_n]$ is a non-decreasing process.
Since both $Y_n(\delta_n)\stackrel{Law}{\rightsquigarrow}0$ and $\limsup_{n\to\infty} \esp{J_t([M_n])}<\infty$ hold,
uniform integrability yields $\limsup_{n\to\infty} E\left\{J_{\delta_n}(Y_n) \right\}=0$ and condition (A) is met by $\langle  M_n  \rangle$.
Proposition \ref{prop:decomp} implies $\esp{J_t(\langle M_n \rangle)}\le \esp{J_t([M_n])}$
so $J_t(\langle  M_n  \rangle)$ is tight for every $t>0$.
Condition (1) of Lemma \ref{lem:aldous} is satisfied as well and $\langle  M_n  \rangle $ is $\cJ_1$-tight.
A similar argument yields the results for $M_n$.

c) See \citet[Theorem II.3.2]{Rebolledo:1979}.

d) See \citet[Corollary VI.6.30]{Jacod/Shiryaev:2003}.
\qed


\section{Proofs of the main results}\label{app:main_results}

\subsection{Proof of Theorem \ref{thm:main_terr1}}\label{pf:thmmain_terr1}

Hypothesis \ref{hyp:An}d implies that $\crochet{ W_n }$ is $\mathcal{C}$-tight, by Proposition \ref{prop:Cincreasing};
applying Theorem \ref{thm:jumps_vanish} to the martingale $W_n = M_n\circ \tau_n$,
taking $B_n$ to be a Brownian motion independent of $M_n$, yields
$(W_n, \crochet{ W_n } ) \stackrel{\mathcal{C}}{\rightsquigarrow} (W,\crochet{ W } )$,
with $W$ a standard Brownian motion since $\crochet{ W }_t=t$. Hypothesis \ref{hyp:An}b.i states that
$A_n \stackrel{\cJ_1}{\rightsquigarrow}A$, so $A_n \stackrel{\cM_1}{\rightsquigarrow} A$ ensues.

In the case where Hypothesis \ref{hyp:An}e.ii or \ref{hyp:An}e.iii holds, so $\tau$ is an inhomogeneous L\'evy process,
$W$ is independent of the pair $(A,\tau)$. The argument runs as follows. Under Hypothesis \ref{hyp:An}e.iii
and by Proposition \ref{prop:inversecontinuity}, the continuity of mapping $A_n\mapsto\tau_n$
under the $\cM_1$-topology on $D_{\uparrow}^0$ ensures that $A_n \stackrel{\cM_1}{\rightsquigarrow} A$ implies both
$\tau_n \stackrel{\cM_1}{\rightsquigarrow} \tau$ and $(A_n,\tau_n) \stackrel{\cM_1}{\rightsquigarrow} (A,\tau)$.
This proof remains unchanged if Hypothesis \ref{hyp:An}e.ii holds instead, starting with
$\tau_n \stackrel{f.d.d.}{\rightsquigarrow} \tau$ --- which implies $\tau_n \stackrel{\cM_1}{\rightsquigarrow} \tau$ ---
and proceeding similarly. Either way, since $(W,\tau)$ is $\mathbb{F}_{\tau}$-adapted, and $\tau$ is a non-decreasing
inhomogeneous L\'evy process, it has no Brownian component, so it is independent of $W$ by Lemma \ref{lem:main}.
Since $A(0)=0$ guarantees the continuity of the reverse mapping $\tau\mapsto A$, $A$ is also independent of $W$
and so is the pair $(A,\tau)$ as the image of Borel measurable mapping $\tau\mapsto(A,\tau)$.
The independence of $W$ and $(A,\tau)$ implies $(A_n,\tau_n,W_n)  \stackrel{f.d.d.}{\rightsquigarrow} (A,\tau,W)$
under either case of Hypothesis \ref{hyp:An}e.ii or \ref{hyp:An}e.iii. Since $W_n \stackrel{\mathcal{C}}{\rightsquigarrow} W$ holds,
$(A_n,\tau_n,W_n) \stackrel{\cM_1}{\rightsquigarrow} (A,\tau,W)$ ensues by Lemma \ref{lem:BasicInequality}.
If instead Hypothesis \ref{hyp:An}e.i holds so $A$ is an inhomogeneous L\'evy process, then the independence of $A$ from $W$
follows from Lemma \ref{lem:main}. Therefore $(A_n,W_n) \stackrel{\cM_1}{\rightsquigarrow} (A,W)$ ensues.
Next, $A_n \stackrel{\cJ_1}{\rightsquigarrow} A$ and $W_n\stackrel{\mathcal{C}}{\rightsquigarrow} W$  together imply
successively $(A_n,W_n) \stackrel{\cJ_1}{\rightsquigarrow} (A,W)$ and
$(W_n\circ A_n,A_n) \stackrel{\cJ_1}{\rightsquigarrow} (W\circ A,A)$, since $(W,\crochet{ W })$ has continuous trajectories,
by the $\cJ_1$-continuity of the composition of functions \citep[Theorem 13.2.2]{Whitt:2002}. It remains to show that $M_n$ converges to $W\circ A$. First, note that under Hypothesis \ref{hyp:An}.a, $0\le\tau_n\circ A_n(t)-t=\sup\{h\ge0; A_n(t+h)=A_n(t)\}<\infty$ holds everywhere,
along almost every trajectory, for any $t\ge 0$. Furthermore,
$$
\tau_n\circ A_n(t)=\sup\{u\ge t; A_n(u)=A_n(t)\}=\inf\{u>t; A_n(u)>A_n(t)\}
$$
is an $\mathbb{F}_{t+\cdot}$-stopping time for every fixed $t\ge0$, since the filtration is right continuous
and $A_n$ is right continuous. Applying Lemma \ref{lem:lenglart} in conjunction with \eqref{eq:next_jump} yields,
for any time $t>0$ and all $\varepsilon,\eta >0$,
$$
P(\sup_{0\le s\le t}|M_n\circ \tau_n\circ A_n(s)-M_n(s)|^2\geq\varepsilon)
\leq \frac{\eta}{\varepsilon} + P(\langle M_n\circ \tau_n\circ A_n \rangle_t - \langle M_n \rangle_t \geq\eta).
$$
Hypothesis \ref{hyp:An}c then yields
$M_n - W_n\circ A_n \stackrel{\mathcal{C}}{\rightsquigarrow} 0$. As a result,
$M_n \stackrel{\cJ_1}{\rightsquigarrow} W\circ A$ and
$(M_n,A_n,W_n) \stackrel{\cJ_1}{\rightsquigarrow} (W\circ A,A,W)$, proving the main statement of the theorem.
Finally, when $J_t(M_n) \stackrel{Law}{\rightsquigarrow} 0$ holds for any $t>0$,
then $M$ is continuous and hence so is $A$, by \citet[Theorem 3.10.2]{Ethier/Kurtz:1986}.
Theorem \ref{thm:jumps_vanish} therefore implies
$(M_n,A_n,W_n) \stackrel{\mathcal{C}}{\rightsquigarrow} (M,A,W)$.
One last appeal to Lemma \ref{lem:BasicInequality} yields
$(M_n,A_n,\tau_n,W_n) \stackrel{\cM_1}{\rightsquigarrow} (M,A,\tau,W)$, completing the proof.
 \qed

\subsection{Proof of Corollary \ref{cor:standardization}}\label{pf:cor:standardization}

The continuity away from the axes, on the complement of
$\disp \cup_{j\in\{1,\ldots,\ell\}}\{(s_1,\cdots,s_{\ell})\in\dR^{\ell}:A_n(s_j)=0\}$, of the map
$$
\{(M_n,A_n)(s_1),\cdots,(M_n,A_n)(s_{\ell})\}\mapsto
\left\{\left(A_n,\frac{M_n}{\sqrt{A_n}}\right)(s_1),\cdots,\left(A_n,\frac{M_n}{\sqrt{A_n}}\right)(s_{\ell})\right\}
$$
from $\dR^{2\ell}$ into itself, for each choice of $0\le s_1 \le \cdots \le s_{\ell} < \infty$,
plus a conditioning argument as in the proof of Equation (\ref{eq:random_rescaling}), together yield the
convergence of the corresponding finite dimensional distributions. The $\mathcal{C}$-continuity of this map when $A$ is
continuous implies the last statement similarly.
\qed

\subsection{Proof of Corollary \ref{cor:main_terr1b}}\label{pf:cor:main_terr1b}

First, build a standard Brownian motion $\hat B$ on a further enlargement of the stochastic basis $(\Omega,\cF,\dF, P)$
so as to be independent of the whole sequence $M_n$. For time change $\disp\epsilon_n(t):=\frac{1-e^{-t}}{n}$,
the perturbed process $\hat M_n:=M_n+\hat B\circ\epsilon_n$ is also a martingale (with respect to the product filtration),
but with both $[\hat M_n]$ and $\langle \hat M_n \rangle$ strictly increasing everywhere.
Writing $\hat\tau_n(s) = \inf\{t\ge0; \langle \hat M_n \rangle _t >s\}$,
note that $\hat\tau_n$ is continuous everywhere and started at $\hat\tau_n(0)=0$.
Furthermore $\disp \sup_{0\le s\le T}|\hat M_n(s)-M_n(s)|\stackrel{\mathcal{C}}{\rightsquigarrow} 0$ holds for every choice of $T>0$.
Put $\hat W_n := \hat M_n\circ \hat \tau_n$. Under the stated conditions of Corollary \ref{cor:main_terr1b},
the martingale $\hat M_n$ also satisfies all Hypothesis \ref{hyp:An}. In fact, Hypotheses \ref{hyp:An}a, \ref{hyp:An}b and \ref{hyp:An}e.i are valid for $\langle \hat M_n \rangle$ with the same $A$. Proposition \ref{prop:inversecontinuity} ensures the validity of
 $\hat\tau_n \stackrel{f.d.d.}{\rightsquigarrow} \tau$ as a consequence of $\langle \hat M_n \rangle \stackrel{f.d.d.}{\rightsquigarrow} A$.
 Hypothesis \ref{hyp:An}c is valid for $(\hat M_n,\crochet{\hat M_n},\hat\tau_n)$,
 irrespective of whether or not it does for $(M_n,A_n,\tau_n)$, since $\crochet{\hat M_n}$ is strictly increasing, which implies
 $\hat\tau_n\circ \crochet{\hat M_n}(t)=t$ for all $t\ge0$. By Lemma \ref{lem:BasicInequality},
$A_n\stackrel{\cJ_1}{\rightsquigarrow} A$ implies $\langle \hat M_n \rangle\stackrel{\cJ_1}{\rightsquigarrow} A$.
Finally, Hypothesis \ref{hyp:An}d ensues from Proposition \ref{prop:Cincreasing}, since the inverse $\tau$ of a strictly increasing $A$
is continuous everywhere, so that $\hat\tau_n \stackrel{\cM_1}{\rightsquigarrow} \tau$ implies successively
$\hat\tau_n \stackrel{\mathcal{C}}{\rightsquigarrow} \tau$ and
$\crochet{\hat W_n}=\langle \hat M_n \rangle\circ\hat\tau_n \stackrel{\cJ_1}{\rightsquigarrow} A\circ\tau$
by the $\cJ_1$-continuity of the composition of functions \citep[Theorem 13.2.2]{Whitt:2002}.
Since $\crochet{\hat M_n}\circ\hat\tau_n(t)\ge t$ is non-decreasing and
$\crochet{\hat W_n}=\crochet{M_n}\circ\hat\tau_n+\epsilon_n^2\circ\hat\tau_n\le \crochet{M_n}\circ\tau_n+\epsilon_n^2\circ\tau_n$,
there ensues $\disp \lim_{n\to\infty}E \crochet{\hat W_n}_t  = \esp{A\circ\tau(t)}$ for any $t\ge 0$
under the uniform integrability of $\crochet{W_n}$.
Therefore, the assumption $\esp{A\circ\tau(t)}=t$ allows for the application of
Theorem \ref{thm:main_terr1} to sequence $(\hat M_n,\langle \hat M_n \rangle,\hat W_n)$,
with eventual limits $M$ and $A$ unchanged, but with a new Brownian motion $\check W$ independent of $A$.
\qed

\subsection{Proof of Theorem \ref{thm:mart_trans}}\label{pf:thm:mart_trans}

Since $(M_n,A_n) \stackrel{\cJ_1}{\rightsquigarrow} (W\circ A,A)$ holds with $W$ a Brownian motion independent of $A$,
there ensues $N_n \stackrel{\cJ_1}{\rightsquigarrow} N_\infty = W\circ A+O_\infty$, by condition a) of Theorem \ref{thm:mart_trans}
and Lemma \ref{lem:BasicInequality}. Applying Corollary \ref{cor:main_terr1b} to $\bar M_n$ yields the existence
of another pair $(\bar W,\bar A_\infty)$ such that
$(\bar M_n,\langle \bar M_n \rangle ) \stackrel{\cJ_1}{\rightsquigarrow} (\bar W\circ \bar A_\infty,\bar A_\infty)$ follows,
provided all the following requirements are met:

\begin{itemize}
\item[(1)] $\langle \bar M_n \rangle _\infty=\infty$ almost surely for each fixed $n\ge1$.
\item[(2)] There is a $D$-valued process $\bar A_\infty$ started at $\bar A_\infty(0)=0$, such that \\
(i) $\langle \bar M_n \rangle \stackrel{\cJ_1}{\rightsquigarrow} \bar A_\infty$ with $\bar A_\infty$ a
strictly increasing inhomogeneous L\'evy process; \\
(ii) $E\langle \bar M_n \rangle _t \to E\bar A_\infty(t)<\infty$ holds for all $t\ge 0$; \\
(iii) $\bar A_\infty(\infty)=\infty$ almost surely.
\item[(3)] $\esp{\bar A_\infty\circ\bar \tau_\infty(t)}=t$ for all $t\ge0$.
\item[(4)] $\langle \bar W_n\rangle_t$ are uniformly integrable for each fixed $t>0$.
\end{itemize}

Item (1) follows from condition b) of Theorem \ref{thm:mart_trans}. Since $\bar A_\infty=\int_0^\cdot \bar Q_\infty^2(u)dA(u)$
is a strictly increasing inhomogeneous L\'evy process by conditions b) and d) of Theorem \ref{thm:mart_trans},
proving Item (2.i) reduces to proving $\langle \bar M_n \rangle \stackrel{\cJ_1}{\rightsquigarrow} \bar A_\infty$.
Since $A_n$ is a $\cJ_1$-tight sequence of non-decreasing processes, it is predictably uniformly tight,
in the sense of \citet[Section VI.6]{Jacod/Shiryaev:2003}. As $\bar Q_n$ is predictable and uniformly bounded,
$\langle \bar M_n \rangle$ is also of bounded variation on compact time sets, under condition b) of Theorem \ref{thm:mart_trans}.
Furthermore, the continuity of $\bar Q_\infty$ in condition c) of Theorem \ref{thm:mart_trans} means that the largest jumps of $\bar Q_n$ and
$\langle \bar M_n \rangle$ go to $0$, since
\[
\sup_{s\in [0,T]}|\Delta \langle \bar M_n \rangle_s|\le A_n(T)\sup_{s\in [0,T]}|\Delta \bar Q_n^2(s)|.
\]
There follows
$\disp \bigl | \langle \bar M_n \rangle - \int_0^\cdot \bar Q_n^2(u-)dA_n(u) \bigr | \stackrel{\mathcal{C}}{\rightsquigarrow} 0$.
By \citet[Theorem VI.6.22]{Jacod/Shiryaev:2003} there ensues
$\disp (\bar Q_n^2, A_n, \langle \bar M_n \rangle ) \stackrel{\cJ_1}{\rightsquigarrow}
\left(\bar Q_\infty^2, A, \bar A_\infty\right)$.

Under condition b) of Theorem \ref{thm:mart_trans}, Item (2.ii) follows by dominated convergence \citep[Proposition App1.2]{Ethier/Kurtz:1986}.
Items (2.iii) and (3) follow from conditions b) and c) of Theorem \ref{thm:mart_trans}, respectively.
Item (4) follows from conditions b) and c) of Theorem \ref{thm:mart_trans}, through
$\disp \bar\tau_n(s) \le\tau_n\left(\frac{s}{\inf_{j,n\ge1}Q_{n,j}^2}\right)$ and therefore
$\disp \langle \bar W_n\rangle_s \le \left(\sup_{j,n\ge1}Q_{n,j}^2\right) A_n\circ\bar\tau_n(s)
\le \left(\sup_{j,n\ge1}Q_{n,j}^2\right) \langle W_n\rangle \circ \left(\frac{s}{\inf_{j,n\ge1}Q_{n,j}^2}\right)$.
Corollary \ref{cor:main_terr1b} yields
$(\bar M_n,\langle \bar M_n \rangle ) \stackrel{\cJ_1}{\rightsquigarrow} (\bar W\circ \bar A_\infty,\bar A_\infty)$.
Finally, $\bar N_n \stackrel{\cJ_1}{\rightsquigarrow} \bar N_\infty=\bar W\circ \bar A_\infty+\bar O_\infty$ ensues from
Lemma \ref{lem:BasicInequality} again, through condition a) of Theorem \ref{thm:mart_trans}; indeed,
since $O_n$ is a $\mathcal{C}$-tight sequence of processes with bounded variation, it is also predictably uniformly tight and the argument
used already to get $\langle \bar M_n \rangle  \stackrel{\cJ_1}{\rightsquigarrow} \bar A_\infty$ also yields
$\bar O_n \stackrel{\mathcal{C}}{\rightsquigarrow} \bar O_\infty$.
\qed

\bibliographystyle{apalike}
\def\cprime{$'$}

\end{document}